\documentclass[12pt]{amsart}
\usepackage{amssymb,amscd}
\usepackage[dvipdfmx]{graphicx,color}

\usepackage{graphicx}

\newtheorem*{Whitney towers}{Theorem~\ref{Whitney towers}}
\newtheorem*{h-towers}{Theorems ~\ref{half} \& \ref{$(n)$-solvable}}

\newtheorem*{surgery curves}{Theorem~\ref{surgery curves}}
\newtheorem*{cg=0}{Theorem~\ref{vanish}}

\newtheorem{thm}{Theorem}[section]

\newtheorem{prop}[thm]{Proposition}
\newtheorem{cla}[thm]{Claim}
\theoremstyle{definition}
\newtheorem{defn}[thm]{Definition}
\newtheorem{note}[thm]{Note}

\newtheorem{prob}[thm]{Problem}

\numberwithin{equation}{section}

\makeatletter
 
 \@addtoreset{figure}{section}
 \makeatother

\newcommand{\x}{\times}

\newcommand{\Z}{\mathbb{Z}}
\newcommand{\N}{\mathbb{N}}

\newcommand{\Q}{\mathbb{Q}}

\newcommand{\R}{\mathbb{R}}

\def\yen{{\setbox0=\hbox{Y}Y\kern-.97\wd0\vbox{hrule height.lex width.98%
\wd0\kern.33ex\hrule height.lex width.98\wd0\kern.45ex}}}

\begin{document}

\title{
Ribbon-move-unknotting-number-two 2-knots, 
pass-move-unknotting-number-two 1-knots,  
and high dimensional analogue}
\author{Eiji Ogasa}


\begin{abstract} 
The (ordinary) unknotting-number of 1-dimensional knots, 
which is defined by using the crossing-change, 
is a very basic and important invariant. 
It is very natural to consider 
the `unknotting-number' 
associated with other local-moves on $n$-dimensional knots ($n\in\N$). 
In this paper we prove the following facts. 
For the ribbon-move on 2-knots, which is a kind of local-move on knots, 
we have the following: 
There is a ribbon-move-unknotting-number-two 2-knot.  
The ribbon-move-unknotting-number of 2-knots is unbounded.       
For the pass-move on 1-knots, which is a kind of local-move on knots, 
we have the following: 
There is a pass-move-unknotting-number-two 1-knot whose (ordinary) unknotting-number  
is $4$.     
For any natural number $n$, there is a 1-knot 
whose pass-move-unknotting-number is$>n$ and whose (ordinary) unknotting-number  
is $4n$.   
For the high-dimensional-pass-move on high-dimensional knots, 
which is a kind of local-move on knots, 
we have the following: 
There is a $(2k+1, 2k+2)$-pass-move-unknotting-number-two  $(4k+2)$-knot.  
The $(2k+1, 2k+2)$-pass-move-unknotting-number of  $(4k+2)$-knot is unbounded.  
There is a $(2k+1, 2k+1)$-pass-move-unknotting-number-two  $(4k+1)$-knot.  
The $(2k+1, 2k+1)$-pass-move-unknotting-number of  $(4k+1)$-knot is unbounded.  
There is a $(4k+1)$-knot whose twist-move-unknotting-number is $n$ 
for any natural number $n$. 
\end{abstract}


\date{}

\maketitle
\section{Introduction}\label{Introduction}




\noindent
The (ordinary) unknotting-number of 1-dimensional knots is a very basic and important invariant of 1-dimensional knots, 
has been studied for a long time, 
and has still many topics to investigate. 
It is well-known 
that there is a 1-knot whose (ordinary) unknotting-number is 
$n$ for any natural number $n$. 
%
The (ordinary) unknotting-number is defined by using the crossing change, which is a local move on knots.  
A local move means as follows: 
When we change a 1-knot $K$ into a 1-knot $J$ by a crossing-change in a 3-ball $B$, 
we make a change only in $B$ and 
that we do not impose any requirement on diffeomorphism type or homeomorphism type 
of $J$ other than the change only in $B$. 
See also Note after Definition \ref{ribbonmove}.

By the way we know other local moves on knots. 
In this paper we discuss the following local-moves: 
\newline
the ribbon-move on 2-dimensional knots, which is defined in \cite{Ogasa04}, 
\newline
the pass-move on 1-knots, which is defined in \cite{Kauffmanon},  
\newline
high-dimensional pass-moves on high-dimensional knots,  which is defined in \cite{Ogasa98n, Ogasa09}, 
and 
\newline
the twist-move on high-dimensional knots, which is defined in \cite{Ogasa09}. 
\newline
(We review their definitions in this paper.) 

It is very natural to ask 
whether there is a knot whose `unknotting-number' associated with 
each of the local moves is two 
and 
whether the `unknotting-number' is unbounded.
(Problems \ref{Arizona}, \ref{Delaware}, \ref{Louisiana}, \ref{Maine}, and \ref{kao}.)
In this paper we give answers to these questions. The each answer is our main theorem. 
Our main results are the following:  \newline
Theorem \ref{Arkansas} about the pass-move on 1-knots,  \newline
Theorem \ref{Florida} about the ribbon-move on 2-knots,   \newline
Theorems \ref{Indiana} and Theorem \ref{Iowa}  
about high-dimensional pass-moves on high-dimensional knots, 
and  \newline
Theorem \ref{kuchi} 
about the twist-move on high-dimensional knots.  \newline
%
%
The statements and the proofs of the first four theorems are different-dimensional analogues of each other. 

$$\text{Table of contents}$$

\hskip1cm
\S\ref{Introduction} {Introduction}

\hskip1cm
\S\ref{California} {The ribbon-move-unknotting-number of 2-knots}

\hskip1cm
\S\ref{Illinois} {The (1,2)-pass-move on 2-knots}

\hskip1cm
\S\ref{Kansas} {Proof of Theorem \ref{Florida}}

\hskip1cm
\S\ref{Kentucky} {Proof of Theorem \ref{Arkansas}}

\hskip1cm
\S\ref{Hawaii} {High-dimensional-pass-moves on high-dimensional knots 

\hskip16.5mm
and their associated `unknotting-number'}

\hskip1cm
\S\ref{Maryland} {Proof of Theorem \ref{Indiana}}

\hskip1cm
\S\ref{Massachusetts} {Proof of Theorem \ref{Iowa}}

\hskip1cm
\S\ref{TM}{The twist-move on high-dimensional knots}

\bigbreak
We review the definitions of the local-moves 
and we state our main theorems. 
We begin by explaining the pass-move on 1-knots 
and the pass-move-unknotting-number of 1-knots.

\bigbreak
We work in the smooth category unless we indicate otherwise. 
Let $n\in\N$. 
If an $n$-(dimensional) oriented submanifold $K\subset S^{n+2}$ is 
 orientation-preserving PL-homeomorphic to the standard sphere $S^n$,  
 $K$ is called an {$n$-$($dimensional$)$ (spherical)-knot}.

Note the following: 
We usually define $n$-knots as above (see e.g. \cite{CochranOrr}). 
Not all $n$-knots are diffeomorphic to the standard $n$-sphere 
although all $n$-knots are PL homeomorphic to the standard $n$-sphere.  
The reason for this is the fact that many exotic $n$-spheres,  
which are not diffeomorphic to the standard $n$-sphere,  
can be embedded smoothly in $S^{n+2}$ 
(see \cite{Levinecob, Levinesimp, Milnor} for the proof of this fact.)

Let $id:S^{n+2}\rightarrow S^{n+2}$ be the identity map. 
We say that $n$-knots $K$ and $K'$ are {\it identical} 
if $id(K)$=$K'$  and 
$id\vert_{K}:K\rightarrow K'$ is an orientation-preserving diffeomorphism map.      
We say that $n$-knots $K$ and $K'$ are {\it equivalent} 
if there exists an orientation-preserving diffeomorphism 
$f:S^{n+2}$ $\rightarrow$ $S^{n+2}$ such that $f(K)$=$K'$  and 
$f\vert_{K}:K\rightarrow K'$ is an orientation-preserving diffeomorphism.    
An $n$-knot $K$ is called a {\it trivial $n$-knot} 
if $K$ is equivalent to the boundary of an $(n+1)$-ball 
embedded in $S^{n+2}$.   


\begin{defn}\label{1pass}
({\bf\cite{Kauffmanon}.}) 
Two 1-knots are {\it pass-move-equivalent}
 if one is obtained from the other 
by a sequence of pass-moves. 
See Figure \ref{koumori}  
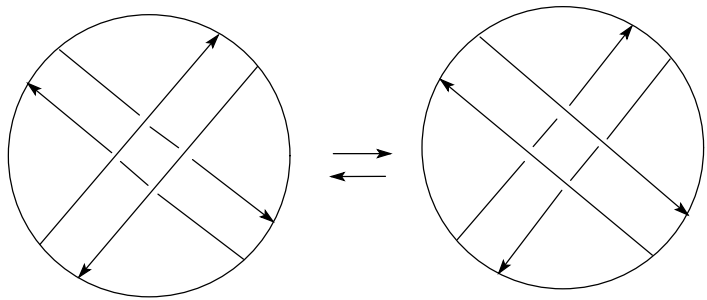
\begin{figure}
\unitlength 0.1in
\begin{picture}(36.39,15.20)(4.01,-22.91)
%
\special{pn 8}%
\special{ar 1138 1552 737 739  0.0000000 6.2831853}%
%
\special{pn 8}%
\special{pa 566 2017}%
\special{pa 1499 921}%
\special{fp}%
\special{sh 1}%
\special{pa 1499 921}%
\special{pa 1441 959}%
\special{pa 1464 962}%
\special{pa 1471 985}%
\special{pa 1499 921}%
\special{fp}%
%
\special{pn 8}%
\special{pa 1705 1085}%
\special{pa 771 2187}%
\special{fp}%
\special{sh 1}%
\special{pa 771 2187}%
\special{pa 829 2149}%
\special{pa 805 2146}%
\special{pa 799 2123}%
\special{pa 771 2187}%
\special{fp}%
%
\special{pn 8}%
\special{pa 1292 1510}%
\special{pa 1143 1401}%
\special{fp}%
%
\special{pn 8}%
\special{pa 1132 1707}%
\special{pa 988 1582}%
\special{fp}%
%
\special{pn 8}%
\special{pa 2375 1660}%
\special{pa 2091 1661}%
\special{fp}%
\special{sh 1}%
\special{pa 2091 1661}%
\special{pa 2158 1681}%
\special{pa 2144 1661}%
\special{pa 2158 1641}%
\special{pa 2091 1661}%
\special{fp}%
%
\special{pn 8}%
\special{pa 2102 1541}%
\special{pa 2401 1541}%
\special{fp}%
\special{sh 1}%
\special{pa 2401 1541}%
\special{pa 2334 1521}%
\special{pa 2348 1541}%
\special{pa 2334 1561}%
\special{pa 2401 1541}%
\special{fp}%
%
\put(18.7000,-21.7000){\makebox(0,0)[lb]{}}%
%
\special{pn 8}%
\special{ar 3303 1510 737 739  0.0000000 6.2831853}%
%
\special{pn 8}%
\special{pa 3148 1504}%
\special{pa 3272 1360}%
\special{fp}%
%
\special{pn 8}%
\special{pa 3298 1733}%
\special{pa 2968 2167}%
\special{fp}%
\special{sh 1}%
\special{pa 2968 2167}%
\special{pa 3024 2126}%
\special{pa 3000 2125}%
\special{pa 2992 2102}%
\special{pa 2968 2167}%
\special{fp}%
%
\special{pn 8}%
\special{pa 3102 1582}%
\special{pa 2751 1991}%
\special{fp}%
%
\special{pn 8}%
\special{pa 3339 1288}%
\special{pa 3664 874}%
\special{fp}%
\special{sh 1}%
\special{pa 3664 874}%
\special{pa 3607 914}%
\special{pa 3631 916}%
\special{pa 3639 939}%
\special{pa 3664 874}%
\special{fp}%
%
\special{pn 8}%
\special{pa 3344 1681}%
\special{pa 3489 1505}%
\special{fp}%
%
\special{pn 8}%
\special{pa 3535 1464}%
\special{pa 3865 1045}%
\special{fp}%
%
\special{pn 8}%
\special{pa 1365 1572}%
\special{pa 1787 1897}%
\special{fp}%
\special{sh 1}%
\special{pa 1787 1897}%
\special{pa 1746 1840}%
\special{pa 1745 1864}%
\special{pa 1722 1872}%
\special{pa 1787 1897}%
\special{fp}%
%
\special{pn 8}%
\special{pa 1087 1339}%
\special{pa 669 997}%
\special{fp}%
%
\special{pn 8}%
\special{pa 1632 2094}%
\special{pa 1184 1753}%
\special{fp}%
%
\special{pn 8}%
\special{pa 942 1536}%
\special{pa 503 1173}%
\special{fp}%
\special{sh 1}%
\special{pa 503 1173}%
\special{pa 542 1231}%
\special{pa 544 1207}%
\special{pa 567 1200}%
\special{pa 503 1173}%
\special{fp}%
%
\special{pn 8}%
\special{pa 2870 930}%
\special{pa 3958 1862}%
\special{fp}%
\special{sh 1}%
\special{pa 3958 1862}%
\special{pa 3920 1803}%
\special{pa 3917 1827}%
\special{pa 3894 1834}%
\special{pa 3958 1862}%
\special{fp}%
%
\special{pn 8}%
\special{pa 3772 2073}%
\special{pa 2664 1153}%
\special{fp}%
\special{sh 1}%
\special{pa 2664 1153}%
\special{pa 2703 1211}%
\special{pa 2705 1187}%
\special{pa 2728 1180}%
\special{pa 2664 1153}%
\special{fp}%
\end{picture}%

\caption{{\bf A pass-move on 1-knots}\label{koumori}}   
\end{figure}
for an illustration of the pass-move. 
 If $K$ and $J$ are pass-move-equivalent and if $K$ and $K'$ are equivalent, 
then we also say that $K'$ and $J$ are pass-move-equivalent. 
\end{defn}

\noindent{\bf Note.} 
\cite{Kauffmanon} proved the following: 
Let $K$ be a 1-knot. 
$K$ is pass-move equivalent to the trivial knot if and only if Arf$K=0$.  
\bigbreak

\begin{defn}\label{1passnum}
Let $K$ be a 1-knot which is pass-move-equivalent to the trivial 1-knot. 
The {\it pass-move-unknotting-number} of $K$ is 
the minimal number of pass-moves which we change $K$ to the trivial 1-knot by. 
\end{defn}

We call the (ordinary) unknotting-number of 1-knots 
the {\it crossing-change-unknotting-number}   
in order to avoid the confusion of notations from now on.

\begin{prop}\label{Alaska}
There is a 1-knot whose pass-move-unknotting-number is one. 
\end{prop}

\begin{figure} 
\vskip-8mm
\hskip-23mm\includegraphics[width=12cm]{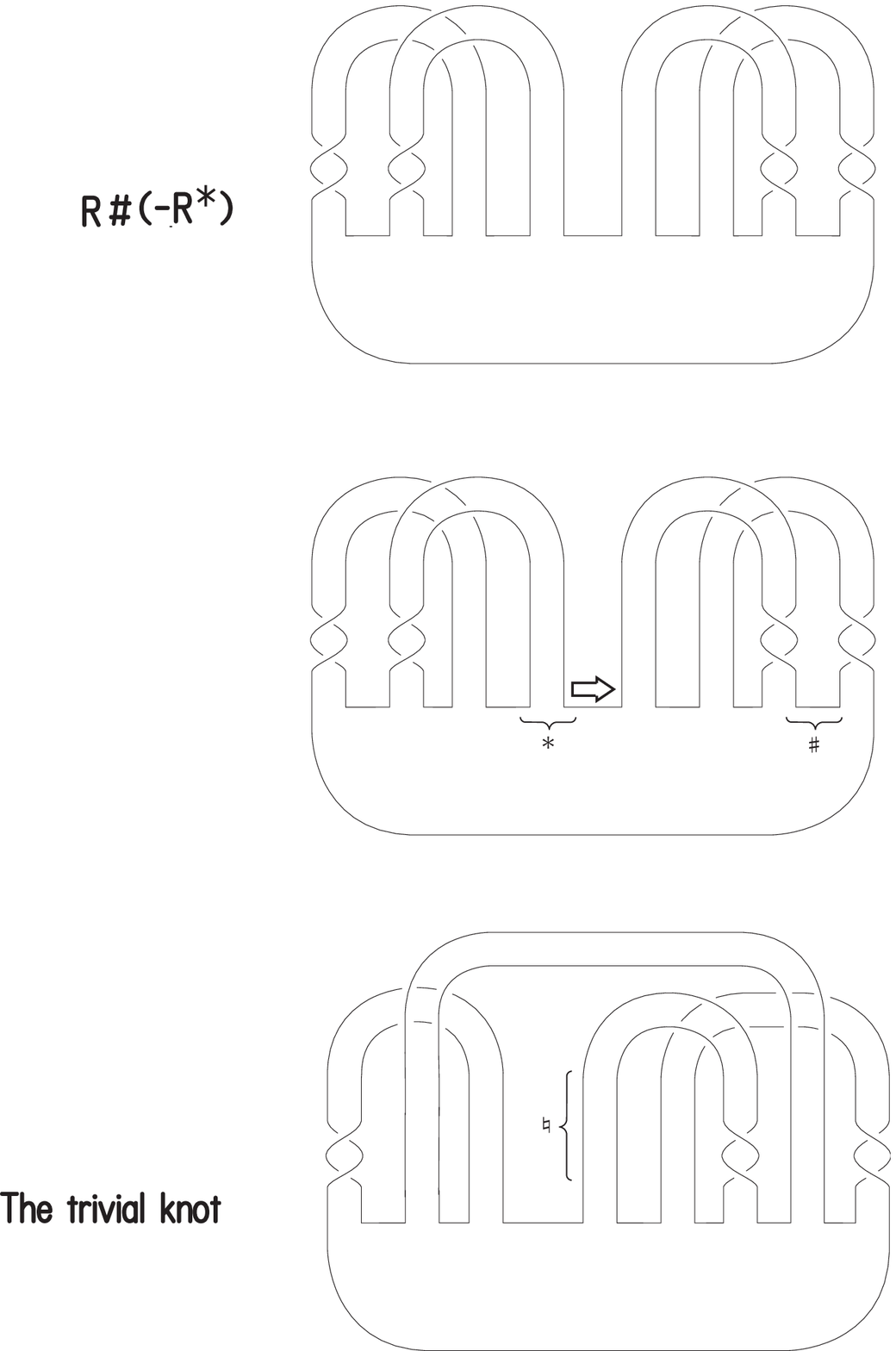}

\caption{{\bf$R\sharp(-R^{*})$ is obtained form the trivial knot \newline
\hskip26mm by a single pass-move.}
\label{kirin}}
\end{figure}

\noindent
{\bf Proof of Proposition \ref{Alaska}.}
Let $R$ be the trefoil knot 
(We do not suppose that $R$ is the right-hand trefoil knot or the left-hand one).    
Then $R\sharp(-R^{*})$ is obtained form the trivial knot by a single pass-move. 
{\it Reason}: 
See Figure \ref{kirin}.  
The  uppermost knot in Figure \ref{kirin}  
is $R\sharp(-R^{*})$. 
See the middle figure in Figure \ref{kirin}.  
We move the part $*$ by isotopy into the direction $\Rightarrow$ to the part $\#$.  
We carry out one pass-move on the resultant 1-knot and 
we obtain the lowermost  knot in Figure \ref{kirin}.  
Note that the pass-move is done near the part $\natural$. 
The readers can check easily by using isotopy  
that the bottom knot is the trivial knot.  

Hence the pass-move-unknotting-number of $R$ is no more than one. 

$R$ is a nontrivial knot because its Alexander polynomial is nontrivial. 
(See Definition \ref{square} for the Alexander polynomial.)
Hence its pass-move-unknotting-number is nonzero. 

Therefore the pass-move-unknotting-number of $R$ is one.  

This completes the proof of Proposition \ref{Alaska}.     \qed\bigbreak

It is very natural to submit the following problem as we state in the first part before 
`Table of contents' of this section. 

\begin{prob}\label{Arizona}   
(1) Is there a pass-move-unknotting-number-two 1-knot 
whose crossing-change-unknotting-number is$\leqq4$?  

\smallbreak\noindent
(2) For any natural number $n$, is there a 1-knot $K$ 
whose pass-move-unknotting-number is$>n$ 
and 
whose crossing-change-unknotting-number is$\leqq4n$?  
\end{prob}

\noindent{\bf Note.} 
It is easy to prove that 
if the crossing-change-unknotting-number of $K$ is$>4n$ and 
the Arf invariant is zero,   
the pass-move-unknotting-number is$>n$. 
Hence we impose the condition on  
the crossing-change-unknotting-number in Problem \ref{Arizona}.     


\bigbreak
We give a positive answer to Problem \ref{Arizona}.(1) (resp. \ref{Arizona}.(2)). 
The answers make one of our main theorems.   

\begin{thm}\label{Arkansas}
$(1)$ There is a pass-move-unknotting-number-two 1-knot 
whose crossing-change-unknotting-number is $4$. 

\smallbreak\noindent
$(2)$ For any natural number $n$, there is a 1-knot 
whose pass-move-unknotting-number is$>n$ 
and 
whose crossing-change-unknotting-number is $4n$. 
\end{thm}



\bigbreak
\section{The ribbon-move-unknotting-number of 2-knots}\label{California}
\noindent  
We use the terms `handle' and `surgeries' in this paper. 
See \cite{Browder, Kirby, Luck, Ranichi, Smale, Wall} 
for the definition of handles (resp. surgeries, the attaching parts of handles, 
the attached part, 
other related terms to handles).  
Note that 
an $a$-dimensional $q$-handle $h^q$ is diffeomorphic to $B^q\x B^{a-q}$ (resp. $B^a$), 
where $B^r$ denotes the $r$-ball, 
and that 
the attaching part of $h^q$ is diffeomorphic $S^{q-1}\x B^{a-q}$. 

\begin{defn}\label{subsur}
Let $x,m\in\N$ and $x<m$. 
Let $X$ be an $x$-dimensional submanifold of an $m$-dimensional manifold $M$.
Suppose that we can embed $X\x [0,1]$ in $M$ so that $X\x\{0\}=X$. 
Suppose that 
an $(x+1)$-dimensional handle $h^p$ is embedded in $M$ and  
is attached to $X\x [0,1]$ ($p\in\N\cup\{0\}, 0\leqq p\leqq x$).  
Suppose that the attaching part of $h^p$ is embedded in $X\x\{1\}$. 
See Figure \ref{bat}.   
\begin{figure} 
\unitlength 0.1in
\begin{picture}(29.70,17.59)(16.60,-24.83)
%
\special{pn 20}%
\special{pa 1660 1834}%
\special{pa 4390 1827}%
\special{pa 4390 2205}%
\special{pa 1667 2205}%
\special{pa 1667 2205}%
\special{pa 1660 1834}%
\special{fp}%
%
\special{pn 8}%
\special{pa 2241 1827}%
\special{pa 2244 1795}%
\special{pa 2247 1763}%
\special{pa 2250 1731}%
\special{pa 2253 1700}%
\special{pa 2256 1668}%
\special{pa 2260 1636}%
\special{pa 2264 1604}%
\special{pa 2268 1572}%
\special{pa 2273 1541}%
\special{pa 2278 1509}%
\special{pa 2284 1477}%
\special{pa 2291 1446}%
\special{pa 2298 1415}%
\special{pa 2307 1384}%
\special{pa 2316 1353}%
\special{pa 2327 1323}%
\special{pa 2338 1293}%
\special{pa 2351 1264}%
\special{pa 2365 1235}%
\special{pa 2381 1208}%
\special{pa 2396 1181}%
\special{pa 2413 1154}%
\special{pa 2429 1127}%
\special{pa 2444 1100}%
\special{pa 2459 1071}%
\special{pa 2473 1042}%
\special{pa 2485 1012}%
\special{pa 2495 979}%
\special{pa 2503 945}%
\special{pa 2511 911}%
\special{pa 2523 882}%
\special{pa 2544 863}%
\special{pa 2574 853}%
\special{pa 2607 844}%
\special{pa 2637 831}%
\special{pa 2665 817}%
\special{pa 2693 802}%
\special{pa 2723 791}%
\special{pa 2753 783}%
\special{pa 2785 777}%
\special{pa 2817 771}%
\special{pa 2848 764}%
\special{pa 2879 757}%
\special{pa 2911 751}%
\special{pa 2942 746}%
\special{pa 2974 739}%
\special{pa 3005 731}%
\special{pa 3036 725}%
\special{pa 3068 724}%
\special{pa 3100 727}%
\special{pa 3132 733}%
\special{pa 3164 740}%
\special{pa 3195 748}%
\special{pa 3225 759}%
\special{pa 3255 771}%
\special{pa 3283 786}%
\special{pa 3310 804}%
\special{pa 3335 824}%
\special{pa 3359 847}%
\special{pa 3380 871}%
\special{pa 3399 897}%
\special{pa 3414 925}%
\special{pa 3427 954}%
\special{pa 3438 985}%
\special{pa 3447 1015}%
\special{pa 3455 1047}%
\special{pa 3463 1078}%
\special{pa 3470 1109}%
\special{pa 3479 1140}%
\special{pa 3488 1171}%
\special{pa 3497 1202}%
\special{pa 3505 1233}%
\special{pa 3514 1264}%
\special{pa 3522 1295}%
\special{pa 3530 1326}%
\special{pa 3537 1357}%
\special{pa 3543 1388}%
\special{pa 3548 1419}%
\special{pa 3553 1450}%
\special{pa 3556 1482}%
\special{pa 3560 1514}%
\special{pa 3563 1546}%
\special{pa 3566 1579}%
\special{pa 3569 1612}%
\special{pa 3572 1644}%
\special{pa 3575 1675}%
\special{pa 3578 1705}%
\special{pa 3580 1736}%
\special{pa 3572 1775}%
\special{pa 3567 1813}%
\special{pa 3585 1827}%
\special{pa 3585 1827}%
\special{sp}%
%
\special{pn 8}%
\special{pa 2535 1834}%
\special{pa 2545 1804}%
\special{pa 2555 1773}%
\special{pa 2564 1743}%
\special{pa 2573 1712}%
\special{pa 2582 1681}%
\special{pa 2591 1650}%
\special{pa 2599 1619}%
\special{pa 2607 1588}%
\special{pa 2614 1557}%
\special{pa 2620 1526}%
\special{pa 2625 1494}%
\special{pa 2630 1462}%
\special{pa 2634 1430}%
\special{pa 2636 1397}%
\special{pa 2640 1365}%
\special{pa 2646 1334}%
\special{pa 2656 1304}%
\special{pa 2670 1275}%
\special{pa 2686 1248}%
\special{pa 2703 1220}%
\special{pa 2719 1193}%
\special{pa 2736 1165}%
\special{pa 2755 1139}%
\special{pa 2778 1116}%
\special{pa 2805 1101}%
\special{pa 2836 1092}%
\special{pa 2868 1086}%
\special{pa 2899 1077}%
\special{pa 2930 1068}%
\special{pa 2960 1058}%
\special{pa 2992 1052}%
\special{pa 3024 1050}%
\special{pa 3056 1053}%
\special{pa 3086 1063}%
\special{pa 3114 1080}%
\special{pa 3138 1101}%
\special{pa 3161 1123}%
\special{pa 3183 1147}%
\special{pa 3202 1172}%
\special{pa 3219 1200}%
\special{pa 3232 1229}%
\special{pa 3243 1260}%
\special{pa 3251 1291}%
\special{pa 3256 1323}%
\special{pa 3259 1355}%
\special{pa 3260 1387}%
\special{pa 3260 1419}%
\special{pa 3261 1451}%
\special{pa 3264 1483}%
\special{pa 3268 1515}%
\special{pa 3274 1546}%
\special{pa 3281 1578}%
\special{pa 3287 1609}%
\special{pa 3291 1641}%
\special{pa 3294 1672}%
\special{pa 3294 1704}%
\special{pa 3292 1736}%
\special{pa 3290 1768}%
\special{pa 3286 1800}%
\special{pa 3284 1820}%
\special{sp}%
%
\special{pn 8}%
\special{pa 1954 2205}%
\special{pa 1958 2237}%
\special{pa 1961 2268}%
\special{pa 1961 2301}%
\special{pa 1956 2333}%
\special{pa 1945 2363}%
\special{pa 1926 2389}%
\special{pa 1908 2415}%
\special{pa 1895 2445}%
\special{pa 1885 2475}%
\special{pa 1884 2478}%
\special{sp -0.045}%
\put(17.3700,-26.5300){\makebox(0,0)[lb]{$X$}}%
%
\special{pn 8}%
\special{pa 2460 1080}%
\special{pa 2428 1082}%
\special{pa 2396 1081}%
\special{pa 2365 1073}%
\special{pa 2335 1061}%
\special{pa 2306 1048}%
\special{pa 2276 1036}%
\special{pa 2246 1024}%
\special{pa 2218 1009}%
\special{pa 2192 990}%
\special{pa 2168 969}%
\special{pa 2160 960}%
\special{sp -0.045}%
\put(19.3000,-9.5000){\makebox(0,0)[lb]{$h^p$}}%
%
\special{pn 8}%
\special{pa 4410 1840}%
\special{pa 4440 1850}%
\special{pa 4464 1871}%
\special{pa 4478 1900}%
\special{pa 4479 1932}%
\special{pa 4481 1964}%
\special{pa 4500 1990}%
\special{pa 4527 2005}%
\special{pa 4561 2010}%
\special{pa 4562 2007}%
\special{pa 4516 2008}%
\special{pa 4497 2032}%
\special{pa 4499 2065}%
\special{pa 4500 2097}%
\special{pa 4490 2125}%
\special{pa 4473 2154}%
\special{pa 4457 2186}%
\special{pa 4432 2200}%
\special{pa 4410 2200}%
\special{sp -0.045}%
\put(46.3000,-20.9000){\makebox(0,0)[lb]{[0,1]}}%
\end{picture}%
\bigbreak\bigbreak\quad\quad
\caption{   
{\bf A handle $h^p$ is attached to $X\x[0,1]$.}   
\label{bat}}
\bigbreak\end{figure}
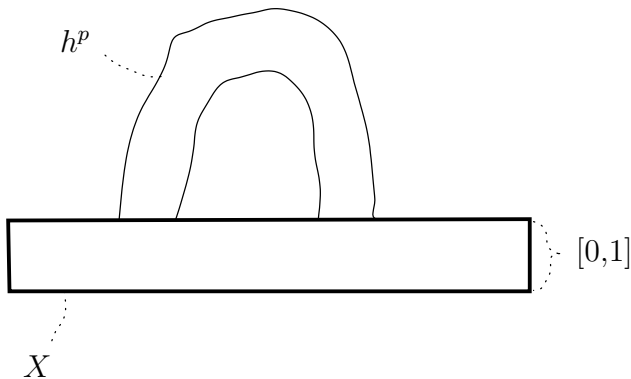
Suppose that $h^p\cap(X\x[0,1])$ is only the attaching part of $h^p$.  
Let 
$X'$=
$\overline{\partial(h^p\cup(X\x[0,1]))-(X\x\{0\})}$.
Note that there are two cases, $\partial X=\phi$ and $\partial X\neq\phi$. 
Then we say that 
 $X'$ is obtained from $X$ 
by {\it the surgery by using the embedded handle $h^p$}. 
We do not say that we use $X\x[0,1]$ if there is no danger of confusion.  
\end{defn}

\noindent
{\bf Note.}  
Of course we can define `embedded surgery' even if we cannot embed $X\x[0,1]$ in $M$. 
However we do not need the case in this paper.  
\bigbreak

\begin{figure}

\includegraphics[width=55mm]{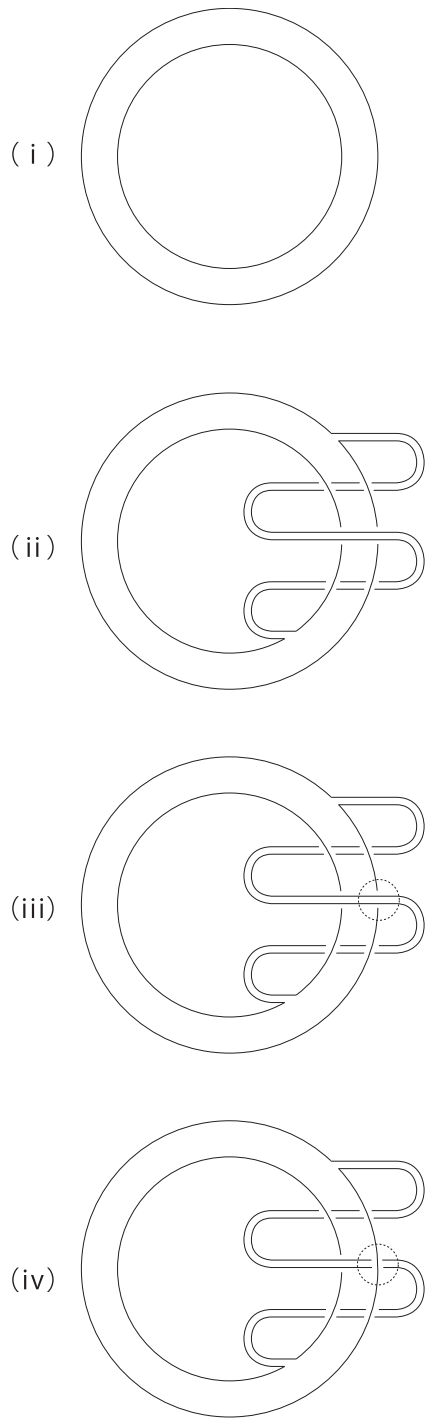}

\vskip-3mm
\caption{ 
 {\bf An example of ribbon-moves on 2-knots } 
\label{spider} }
\end{figure}

We review the definition of the ribbon-move on 2-knots. We begin by showing an example. 
Embed the disjoint union of two copies of $S^2$ 
in $\R^3\x\{t=0\}\subset\R^4\subset S^4$,  
where we regard $\R^4=\R^3\times\{t\in\R\}$, 
as drawn in Figure \ref{spider}.(i).   
Attach an embedded 3-dimensional 1-handle $h^1$ 
to $S^2\amalg S^2$ so that the result of the surgery by using this 1-handle is one $S^2$. 
The 1-handle `rotates' the two $S^2$ as drawn in 
Figure \ref{spider}.(ii).  
If a part of the 1-handle is drawn over (resp. under) a part of the new $S^2$, 
then it means the part of the handle exists in 
$\R^3\x\{t>0\}$ (resp. $\R^3\x\{t<0\}$) as usual. 
The new embedded $S^2$ is embed nontrivially in $\R^4$ because 
the Alexander polynomial is not trivial. 
(See Definition \ref{square} for the Alexander polynomial.)
Note a dotted circle in 
Figure \ref{spider}.(iii),  
which represents the boundary of a 3-ball $B^3$ embedded there. 
We can suppose that (the new $S^2)\cap B^3=$(an annulus)$\amalg$(a disc).  
If we change the over-under of the annulus and the disc in $B^3$, 
then the new $S^2$ becomes a trivial 2-knot as drawn in 
Figure \ref{spider}.(iv).

\begin{defn}\label{ribbonmove}
Let $K_+$ and $K_-$ be 
(not necessarily connected or spherical) 
smooth closed oriented 2-dimensional submanifolds $\subset S^4$. 
%
%
We say that $K_-$ is obtained from $K_+$ by one {\it ribbon-move } 
if there is a 4-ball $B$ embedded 
in $S^4$ with the following properties.

\smallbreak \noindent
(1) 
$K_+$ and $K_-$ differ only in $B$. 

\smallbreak \noindent
(2) 
$B\cap K_+$ (resp. $B\cap K_-$) is diffeomorphic to 
$D^2\amalg (S^1\times [0,1])$, 
where $\amalg$ denotes the disjoint union. 
$B\cap K_+$ (resp. $B\cap K_-$) satisfies the following conditions. 

\begin{figure}
     \begin{center}

    \hskip3cm
\unitlength 0.1in
\begin{picture}(41.74,33.90)(8.50,-34.50)
%
\special{pn 8}%
\special{ar 2829 253 417 186  0.0000000 6.2831853}%
%
\special{pn 20}%
\special{pa 2941 2954}%
\special{pa 2941 2046}%
\special{fp}%
%
\special{pn 20}%
\special{pa 2703 2969}%
\special{pa 2703 2054}%
\special{fp}%
%
\special{pn 20}%
\special{pa 2710 2991}%
\special{pa 2723 2962}%
\special{pa 2750 2946}%
\special{pa 2781 2936}%
\special{pa 2813 2932}%
\special{pa 2845 2934}%
\special{pa 2876 2939}%
\special{pa 2906 2952}%
\special{pa 2928 2974}%
\special{pa 2930 3005}%
\special{pa 2910 3029}%
\special{pa 2881 3042}%
\special{pa 2849 3049}%
\special{pa 2817 3051}%
\special{pa 2786 3048}%
\special{pa 2755 3039}%
\special{pa 2727 3024}%
\special{pa 2711 2996}%
\special{pa 2710 2991}%
\special{sp}%
%
\special{pn 20}%
\special{pa 2703 283}%
\special{pa 2703 1191}%
\special{fp}%
%
\special{pn 20}%
\special{pa 2941 268}%
\special{pa 2941 1183}%
\special{fp}%
%
\special{pn 20}%
\special{ar 2822 246 111 60  0.0000000 6.2831853}%
%
\special{pn 8}%
\special{ar 2829 3051 417 186  0.0000000 6.2831853}%
%
\special{pn 8}%
\special{pa 3253 261}%
\special{pa 3253 3051}%
\special{fp}%
%
\special{pn 8}%
\special{pa 2412 268}%
\special{pa 2412 3051}%
\special{fp}%
%
\special{pn 8}%
\special{ar 1267 261 416 186  0.0000000 6.2831853}%
%
\special{pn 8}%
\special{ar 1267 3058 416 186  0.0000000 6.2831853}%
%
\special{pn 8}%
\special{pa 1691 268}%
\special{pa 1691 3058}%
\special{fp}%
%
\special{pn 8}%
\special{pa 850 276}%
\special{pa 850 3058}%
\special{fp}%
%
\special{pn 8}%
\special{ar 4600 246 416 186  0.0000000 6.2831853}%
%
\special{pn 8}%
\special{ar 4600 3043 416 186  0.0000000 6.2831853}%
%
\special{pn 8}%
\special{pa 5024 253}%
\special{pa 5024 3043}%
\special{fp}%
%
\special{pn 8}%
\special{pa 4183 261}%
\special{pa 4183 3043}%
\special{fp}%
%
\special{pn 20}%
\special{ar 2822 1198 111 60  0.0000000 6.2831853}%
%
\special{pn 20}%
\special{ar 2822 2032 111 59  0.0000000 6.2831853}%
%
\special{pn 20}%
\special{ar 2836 1555 417 186  0.0000000 6.2831853}%
%
\special{pn 8}%
\special{pa 2941 1206}%
\special{pa 4481 1206}%
\special{dt 0.045}%
\special{pa 4481 1206}%
\special{pa 4480 1206}%
\special{dt 0.045}%
%
\special{pn 8}%
\special{pa 2956 2039}%
\special{pa 4496 2039}%
\special{dt 0.045}%
\special{pa 4496 2039}%
\special{pa 4495 2039}%
\special{dt 0.045}%
%
\special{pn 20}%
\special{ar 4570 1206 112 59  0.0000000 6.2831853}%
%
\special{pn 20}%
\special{ar 4570 2046 112 60  0.0000000 6.2831853}%
%
\special{pn 20}%
\special{pa 4451 1213}%
\special{pa 4451 2024}%
\special{fp}%
%
\special{pn 20}%
\special{pa 4696 1236}%
\special{pa 4696 2046}%
\special{fp}%
\put(11.1000,-34.3800){\makebox(0,0)[lb]{t=-0.5}}%
\put(25.9800,-34.3800){\makebox(0,0)[lb]{t=0}}%
\put(43.8400,-34.3800){\makebox(0,0)[lb]{t=0.5}}%
%
\special{pn 8}%
\special{ar 2829 253 417 186  0.0000000 6.2831853}%
%
\special{pn 8}%
\special{pa 2412 268}%
\special{pa 2412 3051}%
\special{fp}%
%
\special{pn 8}%
\special{pa 3253 261}%
\special{pa 3253 3051}%
\special{fp}%
%
\special{pn 8}%
\special{ar 2829 3051 417 186  0.0000000 6.2831853}%
%
\special{pn 20}%
\special{ar 2822 246 111 60  0.0000000 6.2831853}%
%
\special{pn 20}%
\special{pa 2941 268}%
\special{pa 2941 1183}%
\special{fp}%
%
\special{pn 20}%
\special{pa 2703 283}%
\special{pa 2703 1191}%
\special{fp}%
%
\special{pn 20}%
\special{pa 2710 2991}%
\special{pa 2723 2962}%
\special{pa 2750 2946}%
\special{pa 2781 2936}%
\special{pa 2813 2932}%
\special{pa 2845 2934}%
\special{pa 2876 2939}%
\special{pa 2906 2952}%
\special{pa 2928 2974}%
\special{pa 2930 3005}%
\special{pa 2910 3029}%
\special{pa 2881 3042}%
\special{pa 2849 3049}%
\special{pa 2817 3051}%
\special{pa 2786 3048}%
\special{pa 2755 3039}%
\special{pa 2727 3024}%
\special{pa 2711 2996}%
\special{pa 2710 2991}%
\special{sp}%
%
\special{pn 20}%
\special{pa 2703 2969}%
\special{pa 2703 2054}%
\special{fp}%
%
\special{pn 20}%
\special{pa 2941 2954}%
\special{pa 2941 2046}%
\special{fp}%
%
\special{pn 8}%
\special{ar 1267 261 416 186  0.0000000 6.2831853}%
%
\special{pn 8}%
\special{ar 1267 3058 416 186  0.0000000 6.2831853}%
%
\special{pn 8}%
\special{pa 1691 268}%
\special{pa 1691 3058}%
\special{fp}%
%
\special{pn 8}%
\special{pa 850 276}%
\special{pa 850 3058}%
\special{fp}%
%
\special{pn 8}%
\special{ar 4600 246 416 186  0.0000000 6.2831853}%
%
\special{pn 8}%
\special{ar 4600 3043 416 186  0.0000000 6.2831853}%
%
\special{pn 8}%
\special{pa 5024 253}%
\special{pa 5024 3043}%
\special{fp}%
%
\special{pn 8}%
\special{pa 4183 261}%
\special{pa 4183 3043}%
\special{fp}%
\put(23.2000,-36.2000){\makebox(0,0)[lb]{$K_+\cap B$}}%
\put(23.2000,-36.2000){\makebox(0,0)[lb]{$K_+\cap B$}}%
\end{picture}%

\vskip4mm 
F{\footnotesize {IGURE}} \ref{cat}.(1). {\bf Ribbon-move }
\vskip3mm
     \hskip3cm
\unitlength 0.1in
\begin{picture}(41.74,33.70)(8.50,-34.30)
%
\special{pn 8}%
\special{ar 2829 253 417 186  0.0000000 6.2831853}%
%
\special{pn 20}%
\special{pa 2941 2954}%
\special{pa 2941 2046}%
\special{fp}%
%
\special{pn 20}%
\special{pa 2703 2969}%
\special{pa 2703 2054}%
\special{fp}%
%
\special{pn 20}%
\special{pa 2710 2991}%
\special{pa 2723 2962}%
\special{pa 2750 2946}%
\special{pa 2781 2936}%
\special{pa 2813 2932}%
\special{pa 2845 2934}%
\special{pa 2876 2939}%
\special{pa 2906 2952}%
\special{pa 2928 2974}%
\special{pa 2930 3005}%
\special{pa 2910 3029}%
\special{pa 2881 3042}%
\special{pa 2849 3049}%
\special{pa 2817 3051}%
\special{pa 2786 3048}%
\special{pa 2755 3039}%
\special{pa 2727 3024}%
\special{pa 2711 2996}%
\special{pa 2710 2991}%
\special{sp}%
%
\special{pn 20}%
\special{pa 2703 283}%
\special{pa 2703 1191}%
\special{fp}%
%
\special{pn 20}%
\special{pa 2941 268}%
\special{pa 2941 1183}%
\special{fp}%
%
\special{pn 20}%
\special{ar 2822 246 111 60  0.0000000 6.2831853}%
%
\special{pn 8}%
\special{ar 2829 3051 417 186  0.0000000 6.2831853}%
%
\special{pn 8}%
\special{pa 3253 261}%
\special{pa 3253 3051}%
\special{fp}%
%
\special{pn 8}%
\special{pa 2412 268}%
\special{pa 2412 3051}%
\special{fp}%
%
\special{pn 8}%
\special{ar 1267 261 416 186  0.0000000 6.2831853}%
%
\special{pn 8}%
\special{ar 1267 3058 416 186  0.0000000 6.2831853}%
%
\special{pn 8}%
\special{pa 1691 268}%
\special{pa 1691 3058}%
\special{fp}%
%
\special{pn 8}%
\special{pa 850 276}%
\special{pa 850 3058}%
\special{fp}%
%
\special{pn 8}%
\special{ar 4600 246 416 186  0.0000000 6.2831853}%
%
\special{pn 8}%
\special{ar 4600 3043 416 186  0.0000000 6.2831853}%
%
\special{pn 8}%
\special{pa 5024 253}%
\special{pa 5024 3043}%
\special{fp}%
%
\special{pn 8}%
\special{pa 4183 261}%
\special{pa 4183 3043}%
\special{fp}%
%
\special{pn 20}%
\special{ar 2822 1198 111 60  0.0000000 6.2831853}%
%
\special{pn 20}%
\special{ar 2822 2032 111 59  0.0000000 6.2831853}%
%
\special{pn 20}%
\special{ar 2836 1555 417 186  0.0000000 6.2831853}%
\put(11.1000,-34.3800){\makebox(0,0)[lb]{t=-0.5}}%
\put(25.9800,-34.3800){\makebox(0,0)[lb]{t=0}}%
\put(43.8400,-34.3800){\makebox(0,0)[lb]{t=0.5}}%
%
\special{pn 8}%
\special{ar 2829 253 417 186  0.0000000 6.2831853}%
%
\special{pn 8}%
\special{pa 2412 268}%
\special{pa 2412 3051}%
\special{fp}%
%
\special{pn 8}%
\special{pa 3253 261}%
\special{pa 3253 3051}%
\special{fp}%
%
\special{pn 8}%
\special{ar 2829 3051 417 186  0.0000000 6.2831853}%
%
\special{pn 20}%
\special{ar 2822 246 111 60  0.0000000 6.2831853}%
%
\special{pn 20}%
\special{pa 2941 268}%
\special{pa 2941 1183}%
\special{fp}%
%
\special{pn 20}%
\special{pa 2703 283}%
\special{pa 2703 1191}%
\special{fp}%
%
\special{pn 20}%
\special{pa 2710 2991}%
\special{pa 2723 2962}%
\special{pa 2750 2946}%
\special{pa 2781 2936}%
\special{pa 2813 2932}%
\special{pa 2845 2934}%
\special{pa 2876 2939}%
\special{pa 2906 2952}%
\special{pa 2928 2974}%
\special{pa 2930 3005}%
\special{pa 2910 3029}%
\special{pa 2881 3042}%
\special{pa 2849 3049}%
\special{pa 2817 3051}%
\special{pa 2786 3048}%
\special{pa 2755 3039}%
\special{pa 2727 3024}%
\special{pa 2711 2996}%
\special{pa 2710 2991}%
\special{sp}%
%
\special{pn 20}%
\special{pa 2703 2969}%
\special{pa 2703 2054}%
\special{fp}%
%
\special{pn 20}%
\special{pa 2941 2954}%
\special{pa 2941 2046}%
\special{fp}%
%
\special{pn 8}%
\special{ar 1267 261 416 186  0.0000000 6.2831853}%
%
\special{pn 8}%
\special{ar 1267 3058 416 186  0.0000000 6.2831853}%
%
\special{pn 8}%
\special{pa 1691 268}%
\special{pa 1691 3058}%
\special{fp}%
%
\special{pn 8}%
\special{pa 850 276}%
\special{pa 850 3058}%
\special{fp}%
%
\special{pn 8}%
\special{ar 4600 246 416 186  0.0000000 6.2831853}%
%
\special{pn 8}%
\special{ar 4600 3043 416 186  0.0000000 6.2831853}%
%
\special{pn 8}%
\special{pa 5024 253}%
\special{pa 5024 3043}%
\special{fp}%
%
\special{pn 8}%
\special{pa 4183 261}%
\special{pa 4183 3043}%
\special{fp}%
%
\special{pn 8}%
\special{pa 2688 1198}%
\special{pa 1401 1198}%
\special{dt 0.045}%
\special{pa 1401 1198}%
\special{pa 1402 1198}%
\special{dt 0.045}%
%
\special{pn 8}%
\special{pa 2695 2054}%
\special{pa 1423 2054}%
\special{dt 0.045}%
\special{pa 1423 2054}%
\special{pa 1424 2054}%
\special{dt 0.045}%
%
\special{pn 20}%
\special{ar 1274 1198 112 60  0.0000000 6.2831853}%
%
\special{pn 20}%
\special{ar 1274 2039 112 60  0.0000000 6.2831853}%
%
\special{pn 20}%
\special{pa 1155 1206}%
\special{pa 1155 2017}%
\special{fp}%
%
\special{pn 20}%
\special{pa 1401 1228}%
\special{pa 1401 2039}%
\special{fp}%
\put(23.3000,-36.0000){\makebox(0,0)[lb]{$K_-\cap B$}}%
\end{picture}%

\caption{\hskip-2mm(2). \bf  Ribbon-move \label{cat}}
\end{center}
\end{figure}

\smallbreak  
We regard $B$ as 
(a closed 2-disc)$\times[0,1]\times\{t| -1\leqq t\leqq1\}$.
Let $B_t\newline
=$(a closed 2-disc)$\times[0,1]\times\{t \}$.  
Note that $B=\cup B_t$. 
In Figure \ref{cat}.(1) (resp. \ref{cat}.(2)),   
we draw 
$B_{-0.5}$ with $B_{-0.5}\cap K_+$, 
$B_{0}$ with $B_{0}\cap K_+$, 
and 
$B_{0.5}$ with $B_{0.5}\cap K_+$  \newline 
(resp.  
$B_{-0.5}$ with $B_{-0.5}\cap K_-$, 
$B_{0}$ with $B_{0}\cap K_-$, 
and 
$B_{0.5}$ with $B_{0.5}\cap K_-$
).   \newline 
We draw $B_{*}\cap K_+$ and $B_{\#}\cap K_-$ by the bold line, 
where  $*,\#\in\{0.5, 0, -0.5\}$. 
We draw $\partial B_t$ by the fine line.

\smallbreak
$B\cap K_+$ has the following properties:  
$B_t\cap K_+$ is empty for $-1\leqq t<0$ and $0.5<t\leqq1$.
$B_0\cap K_+$ is diffeomorphic to 
$D^2\amalg(S^1\times [0,0.3])\amalg(S^1\times [0.7,1])$. 
$B_{0.5}\cap K_+$ is diffeomorphic to $(S^1\times [0.3,0.7])$. 
$B_t\cap K_+$ is diffeomorphic to $S^1\amalg S^1$ for $0<t<0.5$. 
(Here we draw $S^1\times [0,1]$ to have the corner 
in $B_0$ and in $B_{0.5}$. 
However we can let $B\cap K_+$ in $B$ be a smooth submanifold  
by making the corner smooth naturally.)

\smallbreak
$B\cap K_-$ has the following properties:  
$B_t\cap  K_-$ is empty for $-1\leqq t<-0.5$ and $0<t\leqq1$.
$B_0\cap K_-$ is diffeomorphic to 
$D^2\amalg(S^1\times [0, 0.3])\amalg(S^1\times [0.7, 1])$. 
$B_{-0.5}\cap  K_-$ is diffeomorphic to $(S^1\times [0.3, 0.7])$. 
$B_t\cap  K_-$ is diffeomorphic to $S^1\amalg S^1$ for $-0.5<t<0$. 

\smallbreak
In Figure \ref{cat}.(1) (resp. \ref{cat}.(2))  
there are an oriented cylinder $S^1\times [0,1]$ 
and an oriented disc $D^2$ as we stated above. 
We do not make any assumption about 
the orientation of the cylinder and the disc. 
(Of course it holds that 
this orientation of \newline 
(the cylinder)$\amalg$(the disc) 
coincides with 
the orientation of $B\cap K_+$ (resp. $B\cap K_-$ ).)


\bigbreak
Suppose that $K_-$ is obtained from $K_+$ by one ribbon-move 
and that $K'_-$ is equivalent to $K_-$.   
Then we also say that $K'_-$ is obtained from $K_+$ 
by one {\it ribbon-move}.   
If $K_+$ is obtained from $K_-$ by one ribbon-move,  
then we also say that $K_-$ is obtained from $K_+$ by one {\it ribbon-move}.   
$K_+$ and $K_-$ are said to be {\it ribbon-move equivalent} 
if there are 2-knots \newline
$K_+=\bar{K}_1, \bar{K}_2,...,\bar{K}_{r-1},\bar{K}_r=K_-$, 
where $r$ is a natural number,   
such that 
$\bar{K}_i$ is obtained from $\bar{K}_{i-1}$ $(1< i\leqq r)$ by one ribbon-move. 
\end{defn}

\noindent{\bf Note.}  
When we change a spherical 2-knot $K$ into a closed oriented submanifold $J$ of $S^4$ 
by a ribbon-move in a 4-ball $B$, 
we make a change only in $B$ and 
that we do not impose any requirement on diffeomorphism type or homeomorphism type 
of $J$ other than the change only in $B$.
Note that there are two cases: 
$J$ is diffeomorphic to $S^2$ (resp. $S^2\amalg T^2$). 
This is a reason why we use a term `local' in the term `local-moves' 
as we state in the first part of \S\ref{Introduction}. 
\bigbreak

We explain a derivation of the ribbon-move of 2-knots 
after we review the definition of ribbon 2-knots.
A 2-knot $K$ is called a {\it ribbon} 2-knot  if $L$ satisfies the following properties. 
\smallbreak \noindent
(1) There is a self-transverse immersion 
$f:D^3\rightarrow S^4$ 
such that $f(\partial D^3)=K$.  

\smallbreak \noindent
(2) The singular point set $C$  $(\subset S^4$) 
of $f$ consists of double points. 
$C$ is a disjoint union of 2-discs $D^2_i (i=1,...,k)$. 

\smallbreak \noindent(3) 
Let $j\in\{1,...,k\}$. 
Let $f^{-1}(D^2_j)=D^2_{jB}\amalg D^2_{jS}$. 
The 2-disc $D^2_{jS}$ is 
embedded in the interior of the 3-disc $D^3$.  
The circle $\partial D^2_{jB}$ is 
embedded in 
the boundary of $D^3$.  
The 2-disc $D^2_{jB}$ is embedded in $D^3$.

\smallbreak
It is well-known that it is trivial that  
ribbon 2-knots are changed into the trivial 2-knot by a sequence of ribbon-moves. 
Thus we call the operation defined in Definition \ref{ribbonmove} the ribbon-move.


The author proved the following.

\begin{thm}\label{kyo}   {\bf (\cite{Ogasa04}.)}  
$(1)$ Not all spherical 2-knots are ribbon-move-equivalent to the trivial 2-knot. 

\smallbreak \noindent
$(2)$  There is a nonribbon spherical 2-knot which is ribbon-move-equivalent to the trivial 2-knot. 
\end{thm}

\begin{defn}\label{Connecticut}
Let $K$ be a 2-knot which is ribbon-move-equivalent to the trivial 2-knot. 
The {\it ribbon-move-unknotting-number} of $K$ is 
the minimal number of ribbon-moves which we change $K$ to the trivial 2-knot by. 
\end{defn}

\begin{prop}\label{Georgia}
There is a 2-knot whose pass-move-unknotting-number is one. 
\end{prop}

\begin{figure} 
\hskip-40mm
\includegraphics[width=90mm]{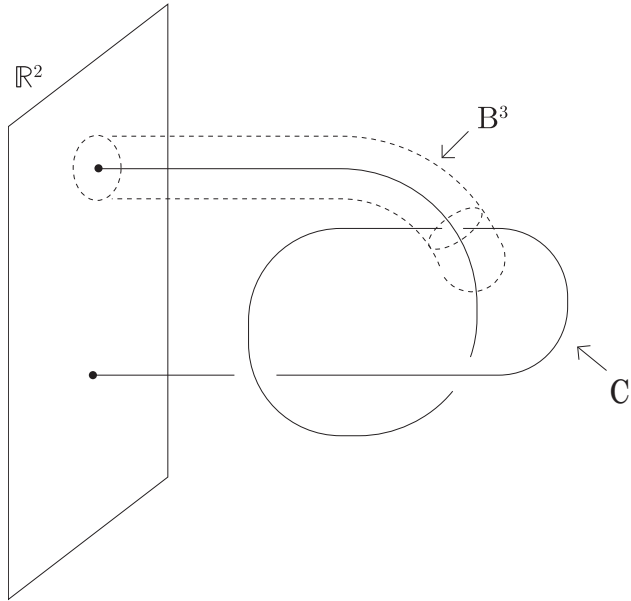}

\caption{  
{\bf An example of spun-knots of 1-knots} \label{hit} }
\end{figure}

\noindent
{\bf Proof of Proposition \ref{Georgia}.}
Figure \ref{spider}.(ii) is an example.  
We give another example. 
It is the spun-knot $S$ of the trefoil knot. 
See Figure \ref{hit}.
%
See \cite{Zeeman} for spun-knots.
This 2-knot $S$ is a nontrivial knot because the Alexander polynomial is nontrivial. 
(See Definition \ref{square} for the Alexander polynomial.)
See the curve $C$ in $R^2\x\{z\geqq0\}$.
Rotate $C$ around $R^2\x\{z=0\}$ as the axis. 
The result is $S$. 
Note that $C\cap (R^2\x\{z=0\})$ consists of two points and is the boundary of $C$. 
Note $B^3$ in $R^2\x\{z\geqq0\}$ which is represented by a dotted curves. 
Note that $B^3\cap C$ is a disjoint union of two curved segments 
and that (two curved segments)$\cap (R^2\x\{z=0\})$ is one point. 
Rotate $B^3$ around $R^2\x\{z=0\}$ as the axis. The result is a 4-ball $B^4$. 
Note that $C\cap B^3$ becomes $S\cap B^4$ 
when we rotate $C$ (resp. $B^3$) around $R^2\x\{z=0\}$.  
Note that $S\cap B^4$ is (a 2-disc)$\amalg$(an annulus). 
We can carry out the ribbon-move in $B^4$. 
This operation changes $S$ to the trivial 2-knot. 
\qed\bigbreak

It is very natural to submit the following problem  
as we state in the first part of \S\ref{Introduction}. 

\begin{prob}\label{Delaware}   
$(1)$ Is there a ribbon-move-unknotting-number-two 2-knot? 

\smallbreak\noindent
$(2)$ For any natural number $n$, is there a 2-knot whose ribbon-unknotting-number is$>n$?  
\end{prob}

We give a positive answer to Problem \ref{Delaware}.(1) (resp. \ref{Delaware}.(2)). 
The answers make one of our main theorems.   

\begin{thm}\label{Florida}
$(1)$ There is a ribbon-move-unknotting-number-two 2-knot. 

\smallbreak\noindent
$(2)$ For any natural number $n$, there is a 2-knot 
whose ribbon-move-unknotting-number is$>n$. 
\end{thm}

\bigbreak
\section{The (1,2)-pass-move on 2-knots}\label{Illinois}
\noindent
In order to prove Theorem \ref{Florida} we review the definition of another local move on 2-knots, 
which is the (1,2)-pass-move on 2-knots defined by the author in \cite{Ogasa04}. 
Why we need the (1,2)-pass-move on 2-knots is because 
we have Proposition \ref{Colorado}, 
which the author proved. 

\begin{defn}\label{12} 
Let $K_+$ and $K_-$ be 2-links in $S^4$. 
We say that 
$K_+$  (resp. $K_-$) is obtained from $K_-$ (resp. $K_+$)  by one {\it $(1,2)$-pass-move } 
if $K_+$ and $K_-$ differ only in a 4-ball $B$ 
embedded in $S^4$ 
with the following properties: 
$B\cap K_+$ is drawn as in Figure \ref{rakuda}.(1).  
\begin{figure}
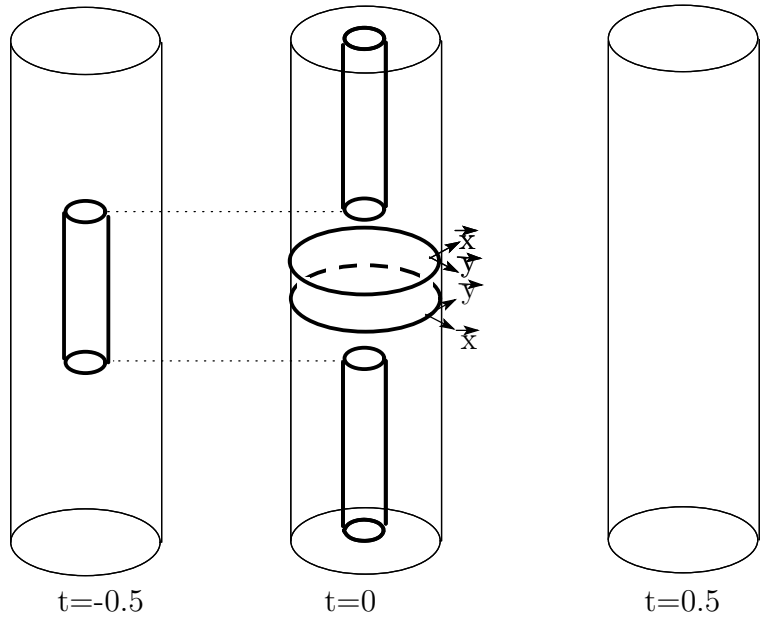

\input 4.1.tex

\bigbreak
\bigbreak
  F{\tiny IGURE} \ref{rakuda}.(1). {\bf The (1,2)-pass-move on 2-knots} 

\bigbreak
\bigbreak
\bigbreak

\input 4.2.tex
\bigbreak
\caption{\hskip-2mm(2). {\bf The (1,2)-pass-move on 2-knots}\label{rakuda}}
\end{figure}
$B\cap K_-$ is drawn as in Figure \ref{rakuda}.(2). 
If $K$ is equivalent to $K'$ and if $K'$ is obtained from $K''$ by a sequence of 
(1,2)-pass-moves, we say that $K$ is {\it $(1,2)$-pass-move-equivalent} to $K''$.  

We draw $B$ as in Definition \ref{ribbonmove}.


$B\cap L_+$ (resp. $B\cap L_-$) is diffeomorphic to 
$D^2\amalg D^2\amalg (S^1\times [0,1])$, 
where $\amalg$ denotes the disjoint union. 
$B\cap L_+$ has the following properties:  
$B_t\cap L_+$ is empty for $-1\leqq t<0$ and $0.5<t\leqq1$.
$B_0\cap L_+$ is  
$(D^2\times\{0.4\})\amalg
(D^2\times\{0.6\})\amalg(S^1\times [0,0.3])\amalg(S^1\times [0.7,1])$. 
$B_{0.5}\cap L_+$ is  $S^1\times [0.3,0.7]$. 
$B_t\cap L_+$ is diffeomorphic to $S^1\amalg S^1$ for $0<t<0.5$. 

$B\cap L_-$ has the following properties:  
$B_t\cap  L_-$ is empty for $-1\leqq t<-0.5$ and $0<t\leqq1$.
$B_0\cap L_-$ is   
$(D^2\times\{0.4\})\amalg
(D^2\times\{0.6\})\amalg(S^1\times [0, 0.3])\amalg(S^1\times [0.7, 1])$. 
$B_{-0.5}\cap  L_-$ is  $S^1\times [0.3, 0.7]$. 
$B_t\cap  L_-$ is diffeomorphic to $S^1\amalg S^1$ for $-0.5<t<0$. 

In Figure \ref{rakuda}.(1) (resp. \ref{rakuda}.(2))   
there are an oriented cylinder $S^1\times [0,1]$ 
and two oriented discs $D^2$. 
We do not make any assumption about 
the orientation of the cylinder.  
We suppose that 
each arrow $\overrightarrow{x}$, $\overrightarrow{y}$ 
in Figure \ref{rakuda}.(1) (resp. \ref{rakuda}.(2)) 
is a tangent vector of each disc at a point. 
(Note we use 
the same notations $\overrightarrow{x}$ (resp. $\overrightarrow{y}$) for different arrows.)
The orientation of each disc in 
Figure \ref{rakuda}.(1) (resp. \ref{rakuda}.(2)) 
is determined by the each ordered set $(\overrightarrow{x},\overrightarrow{y})$. 
The orientation of $B\cap L_+$ (resp. $B\cap L_-$) 
coincides with that of the cylinder and that of the disc. 
\end{defn}

\begin{prop}\label{Colorado}  
{\bf (\cite[Proposition 4.3.(1)]{Ogasa04}.)}
If a 1-knot $K$ is obtained from $J$ by one ribbon-move, 
then $K$ is obtained from $J$ by one $(1,2)$-pass-move.  
\end{prop}


\bigbreak
\section{Proof of Theorem \ref{Florida}}\label{Kansas}
\begin{prop}\label{Mississippi}
Let $K$  
be a 2-knot $\subset S^4$ 
whose ribbon-move-unnknotting-number is one. 
Let $M_3(K)$ be the 3-fold branched covering space of $S^4$ along $K$. 
Then there are three elements$\in H_1(M_3(K);\Z)$ 
which generate $H_1(M_3(K);\Z)$. 
\end{prop}

\begin{figure}
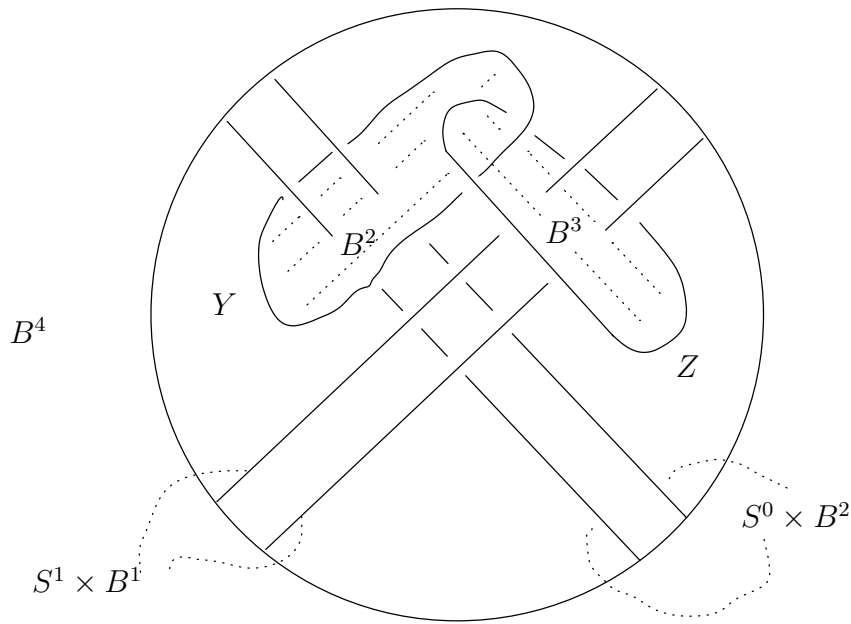
 
\input surgery.tex

\caption{{\bf The (1,2)-pass-move carried out by surgeries}\label{kuma}} 
\end{figure}

\noindent{\bf Proof of Proposition \ref{Mississippi}.}
By Proposition \ref{Colorado} 
 $K$ is obtained from the trivial 2-knot $T$ 
by one (1,2)-pass-move in a 4-ball $B^4\subset S^4$. 
See Figure \ref{kuma}. 
Note that 
$K\cap B^4=(S^0\x B^2)\amalg(S^1\x B^1)$, 
where $\amalg$ denotes the disjoint union. 
Take a 2-ball $B^2$ in the 4-ball $B^4$ such that 
$B^2\cap(S^0\x B^2)$ is two points 
and such that $B^2\cap(S^1\x B^1)=\phi$.  
Call $\partial B^2$, $Y$. 
Take a 3-ball $B^3$ in the 4-ball $B^4$ such that 
$B^3\cap(S^1\x B^1)$ is a circle trivially embedded in $B^3$ 
and that $B^3\cap(S^0\x B^2)=\phi$.  
Call $\partial B^3$, $Z$.  
Suppose that the linking number of $Y$ and $Z$ is one. 
Attach a 5-dimensional 2-(resp.3-)handle to $B^4$ 
along $Y$ (resp. $Z$) with the trivial framing.
Note that these two handles are attached to $S^4$ on time. 
Carry out surgeries by using these two handles on $S^4$. 
Then the new manifold which we obtain is the 4-sphere again, and call it $S^{'4}$.  
Furthermore the new submanifold $\subset S^{'4}$ which is made from $K$  
is the trivial 2-knot $T$. 
 
Note that now we have a compact oriented 5-dimensional manifold $W$ with 
a handle decomposition 

\smallbreak\hskip3mm
$W=(S^4\x[0,1])\cup$(a 5-dimensional 2-handle)

\hskip14mm $\cup$(a 5-dimensional 3-handle)$\cup(S^{'4}\x[0,1])$. 
\smallbreak

Note that $\partial W=(S^4\x\{0\})\amalg(S^{'4}\x\{1\})$.   
Note that there is an embedding map $f:S^2\x [0,1]\hookrightarrow W$ with the following properties: 

\smallbreak
\noindent
(1) 
$f(S^2\x [0,1])\cap(S^4\x\{0\})$ 
is  $f(S^2\x\{0\})$. 
$f(S^2\x [0,1])\cap(S^{'4}\x\{1\})$ 
is  $f(S^2\x\{1\})$. 

$f$ is transverse to $\partial W.$

\smallbreak
\noindent
(2) 
 $f(S^2\x\{0\})$ in $(S^4\x\{0\})$ is $K$. 

$f(S^2\x\{1\})$ in $(S^{'4}\x\{1\})$ is $T$. 

\bigbreak

Take a 3-fold branched covering space $\widetilde{W}$ of $W$ along $f(S^2\x [0,1])$. 
Note the circle which is the core of the attaching part of the 2-handle in 
the above handle decomposition of $W$. 
The circle is null-homologous in $S^4-N(K)$, 
where $N(K)$ is the tubular neighborhood of $K$ in $S^4$.    
Therefore we obtain a compact oriented 5-dimensional manifold $\widetilde{W}$ with 
a handle decomposition 

\smallbreak\hskip3mm
$\widetilde{W}=
(M_3(K)\x[0,1])\cup$(three 5-dimensional 2-handles, $h^2_1$, $h^2_2$, and $h^2_3$)

\hskip14mm 
$\cup$(three 5-dimensional 3-handle, $h^3_1$, $h^3_2$, and $h^3_3$)
$\cup(S^{'4}\x[0,1])$. 
\smallbreak
 
\noindent
Here, 
note that the 3-fold branched covering space of $S^{'4}$ along $T$ is the standard 4-sphere, 
and call it $S^{'4}$ again.

We prove that $\widetilde{W}$ is simply connected. 
{\it Reason.}
Take the dual handle decomposition 

\smallbreak\hskip3mm
$\widetilde{W}=(S^{'4}\x[0,1])\cup$(three 5-dimensional 2-handles, 
$\overline{h^2_1}$, 
$\overline{h^2_2}$, 
$\overline{h^2_3}$)

\hskip14mm
$\cup$(three 5-dimensional 3-handles, 
$\overline{h^3_1}$, 
$\overline{h^3_2}$, 
$\overline{h^3_3}$)
$\cup(M_3(K)\x[0,1])$, 
\smallbreak

\noindent
of the above handle decomposition, where $\overline{h^*_\#}$     
is the dual handle of $h^{5-*}_\#$. 
Take a manifold which is represented by the sub-handle-decomposition  
\smallbreak
$(S^{'4}\x[0,1])\cup$(the three 5-dimensional 2-handles, $h^2_1$, $h^2_2$, and $h^2_3$)
\smallbreak
\noindent
of the dual decomposition of $\widetilde{W}$. 
Since $S^{'4}\x[0,1]$ is simply-connected, this manifold is simply-connected. 
Recall that if we attach 3-handles to a manifold $E$ and we obtain a new manifold $E'$, 
then $\pi_1E\cong \pi_1E'$. 

Therefore $\pi_1(\widetilde{W})=1$. 

Therefore the manifold which is represented by the sub-handle-decomposition 
\smallbreak
$(M_3(K)\x[0,1])\cup$(the three 5-dimensional 2-handles, $h^2_1$, $h^2_2$, and $h^2_3$)
\smallbreak
\noindent 
of the above handle decomposition is simply-connected.

Therefore 
 the cores of the attaching parts of $h^2_1$,  $h^2_2$, and $h^2_3$  
 generate $H_1(M_3(K);\Z)$.     

This completes the proof of Proposition \ref{Mississippi}. \qed\bigbreak

In a similar fashion we can prove the following. 
\begin{prop}\label{gMississippi}
Let $n\in\N$. 
Let $K$ $\subset S^4$ be a 2-knot 
whose ribbon-move-unknotting-number is$\leqq n$. 
Let $M_3(K)$ be the 3-fold branched covering space of $S^4$ along $K$. 
Then there are $3n$ elements$\in H_1(M_3(K);\Z)$ 
which generate $H_1(M_3(K);\Z)$. 
\end{prop}

\begin{defn}\label{Missouri}
Let $n\in\N$. 
Let $K$ be an $n$-knot $\subset S^{n+2}$. 
If $V$ is a connected, compact, oriented, $(n+1)$-dimensional submanifold  $\subset S^{n+2}$ whose boundary is $K$,   
we call $V$ a {\it Seifert hypersurface} for $K$.    
Let $p, n+1-p\in\N$. 
Let $x_1,..., x_\mu$ be $p$-cycles in $V$ 
which compose a basis of $H_p(V;\Z)$/Tor, 
where $\mu\in\N\cup\{0\}$. 
Let $y_1,..., y_\nu$ be $(n+1-p)$-cycles in $V$ 
which compose a basis of $H_{n+1-p}(V;\Z)$/Tor, 
where $\nu\in\N\cup\{0\}$.  
By Poincar\'e duality, we have $\nu=\mu$. 
 Push $y_i$ into the positive (resp. negative) direction of the normal bundle of $V$. 
Call it $y_i^{+}$ (resp. $y_i^{-}$).  
A $(p,n+1-p)$-{\it positive Seifert matrix} for the above submanifold $K$ associated with $V$ represented by 
an ordered basis,  
$(x_1,..., x_\mu)$,   and 
an ordered basis,  
$(y_1,..., y_\mu)$,  
is a $(\mu\x\mu)$-matrix 
$$S=(s_{ij})=({\mathrm{lk}}(x_i, y_j^{+})).$$   
A $(p,n+1-p)$-{\it negative Seifert matrix} 
for the above submanifold $K$ associated with $V$ represented by 
an ordered basis,  
$(x_1,..., x_\mu)$,   and 
an ordered basis,  
$(y_1,..., y_\mu)$,  
is a matrix 
$$N=(n_{ij})=({\mathrm{lk}}(x_i, y_j^{-})).$$   
We have the following: 
Let $S$ and $N$ be as above. 
Then $S-N$ represents the map 
 $\{H_{p}(V;\Z)$/Tor\} $\x \{H_{n+1-p}(V;\Z)$/Tor\} $\rightarrow\Z$ 
which is defined by the intersection product. 
We call  $t\cdot S-N$ the $(p,n+1-p)$-{\it Alexander matrix}   
for $K$ associated with $V$ represented by 
an ordered basis,  
$(x_1,..., x_\mu)$,   and 
an ordered basis,  
$(y_1,..., y_\mu)$. 
`$S$ and $N$' (resp. `$S$ and $t\cdot S-N$',  `$N$ and $t\cdot S-N$')  
are said to be {\it related}  
if 
`$S$ and $N$' (resp. `$S$ and $t\cdot S-N$',  `$N$ and $t\cdot S-N$')  
are defined by using the same $V$, 
the same ordered basis $(x_1,..., x_\mu)$, 
and 
the same ordered basis $(y_1,..., y_\mu)$.   
We sometimes abbreviate 
$(p,n+1-p)$-positive Seifert matrix  
(resp. $(p,n+1-p)$-negative Seifert matrix, $(p,n+1-p)$-Alexander matrix)   
to 
$p$-Seifert matrix  
(resp. $p$-negative Seifert matrix, $p$-Alexander matrix) 
when it is clear from the context. 
\end{defn}

\begin{defn}\label{square}  
Let $n,p\in\N.$
Let $K$ be an $n$-knot$\subset S^{n+2}$. 
Let 
$S_p$ 
$($resp. $N_p)$ 
be  
a $p$-positive 
$($resp. negative$)$  
Seifert matrix for $K$ associated with $V$ represented by 
an ordered basis,  
$(x_1,..., x_\mu)$,   and 
an ordered basis,  
$(y_1,..., y_\mu)$, 
where $\mu\in\N\cup\{0\}$. 
Thus $S_p$ and $N_p$ are related.  

Two polynomials, $f(t)$ and $g(t), \in\Q[t, t^{-1}]$ are said to be 
{\it $\Q[t, t^{-1}]$-balanced } 
if there is an integer $\xi$ 
and a nonzero rational number $r$ 
such that 
$f(t)=r\cdot t^\xi\cdot g(t)$.  

We define 
the {\it $p$-$\Q[t, t^{-1}]$-Alexander polynomial} to be  
the $\Q[t, t^{-1}]$-balanced class of 
`the determinant of $p$-Alexander matrix'  
$${\rm{det}} (t\cdot S_p-N_p).$$ 
\end{defn}

\noindent{\bf Note.} 
This definition is equivalent to 
the spherical-knot-case of Definition 3.1 of \cite{OgasaZ} 
because of Proposition 3.2 of  \cite{OgasaZ}. 



\begin{prop}\label{tsuika}
Let $n,p,n+1-p\in\N.$
Let $N_p$ be 
a $(p,n+1-p)$-{negative Seifert matrix} for $K$ associated with $V$ represented by 
an ordered basis,  
$(x_1,..., x_\mu)$,   and 
an ordered basis,  
$(y_1,..., y_\mu)$, 
where $\mu\in\N\cup\{0\}$.    
Let $S_{n+1-p}$ be 
a $(n+1-p,p)$-{positive Seifert matrix} for $K$ associated with $V$ represented by 
an ordered basis,  
$(y_1,..., y_\mu)$,   and 
an ordered basis,  
$(x_1,..., x_\mu)$. Then we have 
$$N_p=(-1)^{p\cdot n+1}\cdot S_{n+1-p}.$$ 
\end{prop}

\noindent
{\bf Proof of Proposition \ref{tsuika}.}
By the definition of 
$x_i^{+}$ and  $y_i^{-}$, 
${\mathrm{lk}}(y_i, x_j^{+})$ 
$={\mathrm{lk}}(y_i^{-}, x_j)$. 
By \cite[page 541]{Levinepol}, 
${\mathrm{lk}}(y_i^{-}, x_j)$
$=(-1)^{p(n+1-p)+1}{\mathrm{lk}}(x_j, y_i^{-})$. 
Note that $p(1-p)$ is an even number.\qed

\bigbreak
Proposition \ref{tsuika} implies Proposition \ref{Mars}. 

\begin{prop}\label{Mars} 
Let $m\in\N\cup\{0\}.$  
Let $K$ be a $(2m+1)$-dimensional closed oriented submanifold  $\subset S^{2m+3}$. 
Let $S$ be an $(m+1,m+1)$-Seifert matrix. 
Then we have 
$$S=(-1)^{m}\cdot ^t \hskip-1mm S.$$ 
\end{prop}


\bigbreak
We use  the following proposition, which the Mayer-Vietoris sequence implies. 

\begin{prop}\label{Alabama}  
{\bf(Folklore.)} 
Let $n\in\N.$  
Let $A$ be an $n$-knot. 
Let an $(l\times l)$-matrix $Z$ be a positive-$p$-Seifert matrix for $A$, 
where $l\in\N\cup\{0\}$ and $p\in\N$. 
Suppose that $Z$ is invertible as $\Z$-valued matrix. 
Let $k\in\N$. 
Let $X_k(A)$ be the $k$-fold branched covering space of $S^{n+2}$ along $A$.
Then $H_p(N_k(A);\Z)$ is generated by 
$(^t\hskip-1mm Z Z^{-1})^k-I$, where $I$ is the $(l\times l)$-identity matrix. 
%
\end{prop}


Let $R_i$ be the trefoil knot for $i=1,2$ 
(We do not suppose 
whether $R_i$ is the right-hand trefoil knot or the left-hand one for each $i$, nor 
whether $R_1$ is equivalent to $R_2$).      
Let $P$ be a spun-knot of  $R_1\sharp R_2$.  
Note that $P$ is a 2-knot$\subset S^4$. 
It is well-known that spun-knots are ribbon-knots. 
Hence $P$ is ribbon-move equivalent to the trivial 2-knot.

There is a Seifert surface $V$ for $R_i$ ($i=1,2$) with the following properties: 

\smallbreak\noindent (1)  
$H_1(V;\Z)\cong\Z\oplus\Z$. 

\smallbreak\noindent (2)    
There is  an ordered set $(x_1,x_2)$ of basis of $H_1(V;\Z)\cong\Z\oplus\Z$. 
The intersection matrix $(x_k\cdot x_l)$ ($k, l\in\{1,2\}$)  
on $H_1(V;\Z)$ is 
$
\begin{pmatrix}
0&1\\
-1&0
\end{pmatrix}. 
$
Note that Poincar\'e dual of $x_1$ is $x_2$.

\smallbreak\noindent (3)  
The Seifert matrix $({\rm lk}(x_k, x_l^+))$   for $R_i$ ($i=1,2$) is    
$X=
\begin{pmatrix}
-1&1\\
0&-1
\end{pmatrix}
$

\smallbreak 

Therefore we have the following: 

\smallbreak
$
X^{-1}=
\begin{pmatrix}
-1&-1\\
0&-1
\end{pmatrix}, 
$
$
^t\hskip-1mm X X^{-1}=
\begin{pmatrix}
1&1\\
-1&0
\end{pmatrix}, 
$
$
(^t\hskip-1mm X X^{-1})^2=
\begin{pmatrix}
0&1\\
-1&-1
\end{pmatrix}, 
$
$
(^t\hskip-1mm X X^{-1})^3=
\begin{pmatrix}
-1&0\\
0&-1
\end{pmatrix}.  
$
%
%

\bigbreak
By the definition of the spun-knot, 
$P$ has a Seifert hypersurface $V$ as follows: 

\smallbreak\noindent (1)  
$H_1(V;\Z)\cong\Z\oplus\Z$.   
$H_2(V;\Z)\cong\Z\oplus\Z$. 

\smallbreak\noindent (2)    
There is an ordered set  $(x_1,x_2)$ of basis of $H_1(V;\Z)\cong\Z\oplus\Z$. 
There is an ordered set  $(y_1,y_2)$ of basis of $H_2(V;\Z)\cong\Z\oplus\Z$. 
The intersection matrix $(x_k\cdot y_l)$ ($k, l\in\{1,2\}$)  
on $H_1(V;\Z)$ (resp. $H_2(V;\Z)$) is 
$
\begin{pmatrix}
0&1\\
-1&0
\end{pmatrix}. 
$
Note that Poincar\'e dual of $x_1$ (resp. $x_2$) is $y_2$ (resp. $-y_1$).

\smallbreak\noindent (3)  
The Seifert matrix $({\rm lk}(x_k, y_l^+))$   for $R_i$ ($i=1,2$) is  
$
\begin{pmatrix}
-1&1\\
0&-1
\end{pmatrix}. 
$

\bigbreak


Let $M_3(P)$ be the 3-fold branched covering space of $S^4$ along $P$. 

By Proposition \ref{Alabama} 
we have $H_1(M_3(P);\Z)\cong\Z_2\oplus\Z_2\oplus\Z_2\oplus\Z_2$. 
Hence we need no less than four generators in order to generate 
$H_1(M_3(P);\Z)$. 

Suppose that the ribbon-move-unknotting-number of $P$ is$\leqq1$. 
By Proposition \ref{Mississippi},  
we can take three generators in order to generate $H_1(M_3(P);\Z)$. 
We arrived at a contradiction. 
Hence the ribbon-move-unknotting-number of $P$ is$\geqq2$. 

The ribbon-move-unknotting-number of $P$ is$\leqq2$.
{\it Reason}: See Proof of Proposition \ref{Georgia}.  

Hence the ribbon-move-unknotting-number of $P$ is two.

This completes the proof of Theorem \ref{Florida}.(1). 
\bigbreak

Let $n\in\N$. 
Let $m\in\N$ and $\frac{2m}{3}>n$.  
Let $\#^mP$ be the connected-sum of $m$-copies of $P$. 
Since $P$ is ribbon-move equivalent to the trivial 2-knot, 
 $\#^mP$ is ribbon-move equivalent to the trivial 2-knot. 

Let $N_3(\#^mP)$ be the 3-fold branched covering space of $S^4$ along $\#^mP$. 
By Proposition \ref{Alabama} 
we have $H_1(M_3(\#^mP);\Z)\cong\oplus^{2m}\Z_2$. 
Hence 
we need no less than $2m$ generators in order to generate $H_1(M_3(P);\Z)$.

Suppose that  the ribbon-move-unknotting-number of $\#^mP$ is $\leqq n$.
By Proposition \ref{gMississippi} 
$H_1(M_3(\#^mP);\Z)$ can take $3n$ generators. 
Since $2m>3n$, we arrived at a contradiction. 
Therefore  the ribbon-move-unknotting-number of $\#^mP$ is$>n$. 

This completes the proof of Theorem \ref{Florida}.(2). 

This completes the proof of Theorem \ref{Florida}.\qed

\bigbreak
\section{Proof of Theorem \ref{Arkansas}}\label{Kentucky}
 
%
%

\begin{prop}\label{six}  
Let $J$  
be a 1-knot $\subset S^3$ 
whose pass-move-unknotting-number is one. 
Let $N_3(J)$ be a 3-fold branched covering space of $S^3$ along $J$. 
Then there are six elements$\in H_1(N_3(J);\Z)$ 
which generate $H_1(N_3(J);\Z)$. 
\end{prop}

\noindent {\bf Proof of Proposition \ref{six}.}   
Take a 3-ball $B^3\subset S^3$ 
where we carry out the pass-move which changes $J$ into $T$. 
See Figure \ref{tako}.  
\begin{figure} 
\vskip-11mm
\input pass.tex

\smallbreak
\caption{{\bf The pass-move carried out by surgeries}\label{tako}}
\smallbreak
\end{figure}
Note that $J\cap B^3$ is regarded as \newline $(S^0\x B^1)\amalg(S^0\x B^1)$. 
Call one of the two $S^0\x B^1$, $A_1$, and the other $A_2$. 
Take two 2-balls, $B^2_1$ and $B^2_2$, in the 3-ball $B^3$ 
such that $B^2_i\cap A_i$ is two points 
and such that $B^2_i\cap A_j=\phi$ 
($i=1,2$, and $i\neq j$). 
Call $\partial B^2_i$, $Y_i$ ($i=1,2$).   
Suppose that the linking number of $Y_1$ and $Y_2$ is one. 
Attach a 4-dimensional 2-handle to $B^3$ 
along $Y_i$ with the trivial framing ($i=1,2$).  
Note that these two handles are attached to $S^3$ on time. 
Carry out surgeries by using these two handles on $S^3$. 
Then the new manifold which we obtain is the 3-sphere again, and call it $S^{'3}$.  
Furthermore the new submanifold$\subset S^{'3}$ which is made from $J$ 
is the trivial 1-knot $T$. 
 
Note that we now have a compact oriented 4-dimensional manifold $U$ with 
a handle decomposition 

\smallbreak\hskip3mm
$U=(S^3\x[0,1])\cup$(two 3-dimensional 2-handles)
$\cup(S^{'3}\x[0,1])$. 
\smallbreak

Note that $\partial U=(S^3\x\{0\})\amalg(S^{'3}\x\{1\})$. 
There is an embedding map 
$f:S^1\x [0,1]\hookrightarrow U$ with the following properties: 

\smallbreak
\noindent
(1) 
$f(S^1\x [0,1])\cap(S^3\x\{0\})$ 
is  $f(S^1\x\{0\})$. 
$f(S^1\x [0,1])\cap(S^{'3}\x\{1\})$ 
is  $f(S^1\x\{1\})$. 

$f$ is transverse to $\partial U.$

\smallbreak
\noindent
(2) 
$f(S^1\x\{0\})$ in $S^3\x\{0\}$ is $J$. 

$f(S^1\x\{1\})$ in $S^{'3}\x\{1\}$ is $T$. 

\bigbreak

Take a 3-fold branched covering space $\widetilde{U}$ of $U$ along $f(S^1\x [0,1])$. 
Note the circle which is the core of the attaching part of each of the two 2-handles in 
the above handle decomposition of $U$. 
Each of the two circles is null-homologous in $S^3-N(J)$, 
where $N(J)$ is the tubular neighborhood of $J$ in $S^3$.    
Therefore we obtain a compact oriented 4-dimensional manifold $\widetilde{U}$ with 
a handle decomposition 

\smallbreak\hskip3mm
$\widetilde{U}=(N_3(J)\x[0,1])\cup$
(six 4-dimensional 2-handles, $h^2_1$,...,$h^2_6$)
$\cup(S^{'3}\x[0,1])$.  
\smallbreak
 
\noindent 
Here, note that a 3-fold branched covering space of $S^3$ along $T$ is the standard 3-sphere, 
and call it $S^{'3}$ again.

We prove that $\widetilde{U}$ is simply connected. 
{\it Reason.}
Take the dual handle decomposition 

\smallbreak\hskip3mm
$\widetilde{U}=(S^{'3}\x[0,1])\cup$
(six 4-dimensional 2-handles, $\overline{h^2_1}$,,...,$\overline{h^2_6}$)$\cup(N_3(J)\x[0,1])$, 
\smallbreak

\noindent
of the above handle decomposition, where $\overline{h^2_\#}$ 
is the dual handle of $h^2_\#$. 
Since $S^{'3}\x[0,1]$ is simply-connected, 
$\widetilde{U}$ is simply-connected.




Therefore 
 the cores of the attaching parts of $h^2_1$,...,$h^2_6$ generate  
   $H_1(N_3(J);\Z)$.   

This completes the proof of Proposition \ref{six}.  \qed  \bigbreak 

In a similar way, we can prove the following.

\begin{prop}\label{gsix}
Let $n\in\N$. 
Let $J$ $\subset S^3$ be a 1-knot 
whose pass-move-unknotting-number is$\leqq n$.   
Let $N_3(J)$ be a 3-fold branched covering space of $S^3$ along $J$. 
Then there are 6$n$ elements$\in H_1(N_3(J);\Z)$ 
which generate $H_1(N_3(J);\Z)$. 
\end{prop}

Let $R$ be the trefoil knot 
(We do not suppose that $R$ is the right-hand trefoil knot or the left-hand one).    
Let $C=(R\# (-R^*))\#(R\# (-R^*))$. 
Note that Arf $C=0$. 
By Note to Definition \ref{1pass}, $C$ is pass-move equivalent to the trivial 1-knot.

By Proposition \ref{Alabama} and the calculations right after Proposition \ref{Alabama},   
$H_1(N_3(C);\Z)\cong\oplus^8\Z_2$. 
Hence we need no less than eight generators to generate $H_1(N_3(C);\Z)$. 

Suppose that 
the pass-move-unknotting-number of $C$ is$\leqq1$. 
By Proposition \ref{six} 
$H_1(N_3(C);\Z)$ can take six generators. 
We arrived at a contradiction. 

Therefore the pass-move-unknotting-number of $C$ is$\geqq2$. 

The pass-move-unknotting-number of $C$ is$\leqq2$. 
{\it Reason}: See Proof of Proposition \ref{Alaska}. 

Therefore the pass-move-unknotting-number of $C$ is two. 

It is trivial to prove that the crossing-change-unknotting-number of $R$ is one. 

The crossing-change-unknotting-number of $C$ is $4$ because of 
\cite[Proof of Theorem 10.1  in page 420 and (2.4) in page 389]{Murasugi}.

This completes the proof of Theorem \ref{Arkansas}.(1). 

\bigbreak

Let $n\in\N$.   
Let $\#^nC$ be the connected-sum of $n$ copies of $C$. 
Since $C$ is pass-move equivalent to the trivial 1-knot, 
 $\#^nC$ is pass-move equivalent to the trivial 1-knot. 

Let $N_3(\#^nC)$ be the 3-fold branched cyclic covering space of $S^3$ along $\#^nC$.
By Proposition \ref{Alabama}, 
$H_1(N_3(\#^nC);\Z)\cong\oplus^{8n}\Z_2$. 
Hence we need no less than $8n$ generators which generate  $H_{2k+1}(N_3(\#^nC);\Z)$. 

Suppose that the pass-move-unknotting-number of $\#^nC$ is$\leqq n$. 
By Proposition \ref{gsix}, we can prove that  
$H_1(N_3(\#^nC);\Z)$ can take $6n$ generators.  
We arrived at a contradiction. 

Therefore the pass-move-unknotting-number of $\#^nC$ is$>n$. 


The crossing-change-unknotting-number of $\#^nC$ is $4n$ because of   
\cite[Proof of Theorem 10.1  in page 420 and (2.4) in page 389]{Murasugi}.

This completes the proof of Theorem \ref{Arkansas}.(2). 

This completes the proof of Theorem \ref{Arkansas}. \qed 


\bigbreak
\section
{High-dimensional-pass-moves on high-dimensional knots 
and their associated `unknotting-number'}\label{Hawaii}

\noindent
Local moves on high dimensional knots were defined in  
\cite{Ogasa98n,  Ogasa04, Ogasa09}.  
They have been  researched in  
\cite{KauffmanOgasa,  KauffmanOgasaII,  KauffmanOgasaB, 
Ogasa98n,  Ogasa02,  Ogasa04, Ogasa07,   Ogasa09, OgasaT3, OgasaIH}. 
We show an example of them before 
we review the definition of high-dimensional pass-moves on high dimensional knots.

\bigbreak
\noindent{\bf Lemma.} {\it  
Let $p\in\N$. 
Letting $B^p$ denote a $p$-dimensional ball, 
we can write 
$$S^p=B^p_u\cup B^p_d$$  
\hskip4cm$S^p\x S^q
=(B^p_u\cup B^p_d)\x(B^q_u\cup B^q_d).$

Thus 
$$S^p\x S^q=(B^p_u\x B^q_u)\cup(B^p_u\x B^q_d)\cup(B^p_d\x B^q_u)\cup(B^p_d\x B^q_d).$$
}
\noindent{\bf Proof.}  
Use the fact  
\vskip3mm
\hskip45mm$(X\cup Y)\times Z=(X\times Z)\cup(Y\times Z).$  \qed

\smallbreak 
Let $p,q\in\N$. 
Let 
$$F= (S^p\x S^q)-\text{Int}(B^p_u\x B^q_u).$$  
We indicate $F$ in the figure below and abbreviate $B^\sharp_\star$ to  $B_\star$. 

\bigbreak 
\includegraphics[width=12cm]{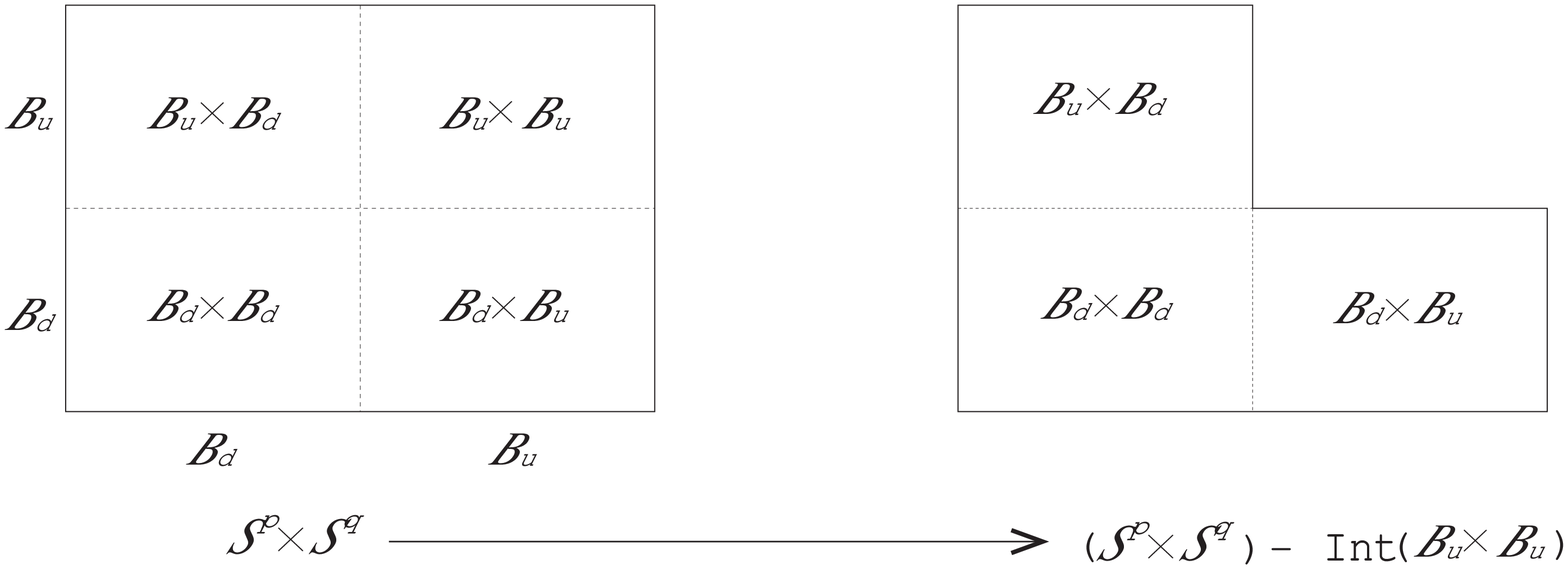}  
\bigbreak

$F$ is drawn 
in another way as below. 
Note that we can bend the corner of $B^p_u\x B^q_u$ 
and change it into the $(p+q)$-dimensional ball. 
Let $p+q=n+1$. 
Hence the boundary of $F$ is $S^n$.

\bigbreak
\includegraphics[width=14cm]{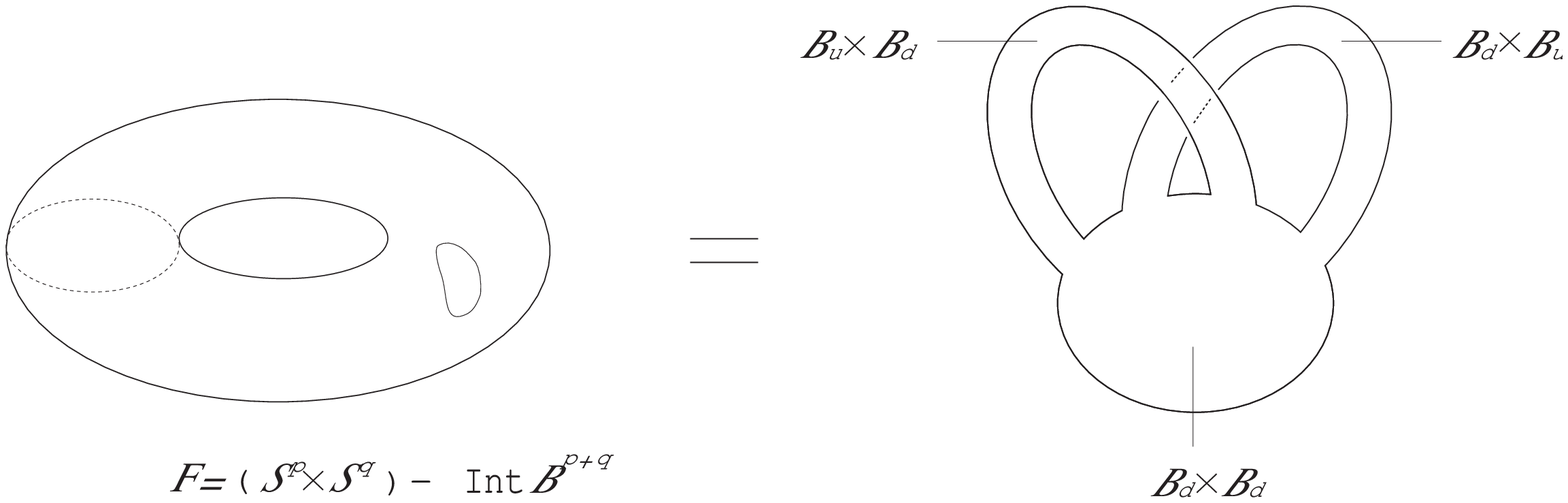}  
\smallbreak
\hskip4cm F{\tiny IGURE} \ref{Hawaii}.1:{\bf $(S^p\x S^q)-$Int$B^{p+q}$} 
\bigbreak

We can regard $B^p_d\x B^q_d$ as a $(p+q)$-dimensional 0-handle, \newline
$B^p_u\x B^q_d$ as a $(p+q)$-dimensional $p$-handle, and  \newline
$B^p_d\x B^q_u$ as a $(p+q)$-dimensional $q$-handle.

\noindent
Embed $F\subset S^{n+2}$ as follows: 
%
%
Embed $B^{p+1}\subset S^{n+2}$. 
Take the tubular neighborhood $N(\partial B^{p+1})$ of $\partial B^{p+1}$ in $S^{n+2}$. 
Take $\partial(N(\partial B^{p+1}))$. 
Note that $\partial(N(\partial B^{p+1}))$ is diffeomorphic to $S^p\x S^q$.
Embed $F\subset S^p\x S^q$ as above. 

\bigbreak
\hskip3cm\includegraphics[width=5cm]{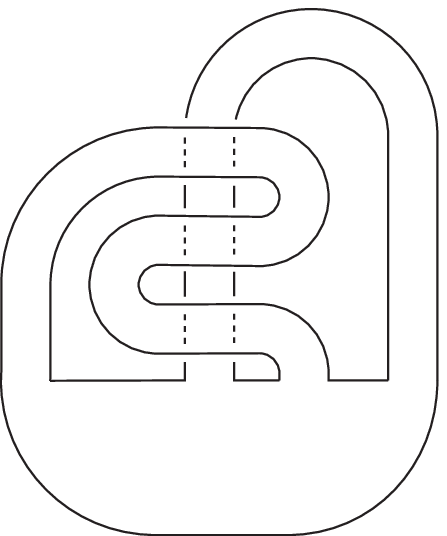}   
\smallbreak
\hskip4cm F{\tiny IGURE} \ref{Hawaii}.2:  {\bf A trivial $n$-knot} 

\bigbreak

The boundary of $F$ in  $S^{n+2}$ is an $n$-knot.  
Furthermore it is the trivial $n$-knot. 
Carry out a `local move' on this $n$-knot 
in an $(n+2)$-ball, which is denoted by a dotted circle  
in the following figure.

\bigbreak
\hskip4cm\includegraphics[width=5cm]{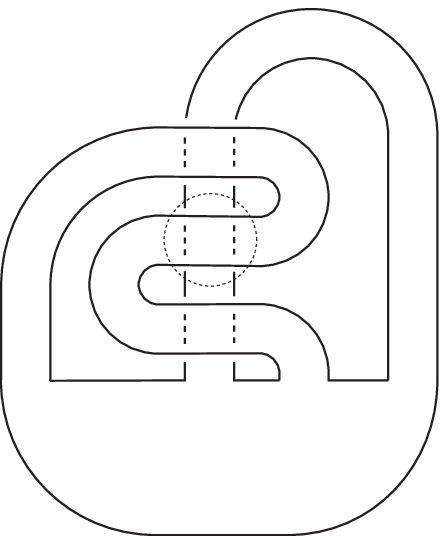}  
\smallbreak
\hskip2cm   F{\tiny IGURE} \ref{Hawaii}.3:{\bf A local move will be carried out in the dotted  

\hskip20mm $(n+2)$-ball. The resulting $n$-knot $K$ is a nontrivial $n$-knot.}  
\bigbreak


\bigbreak
\hskip4cm\includegraphics[width=5cm]{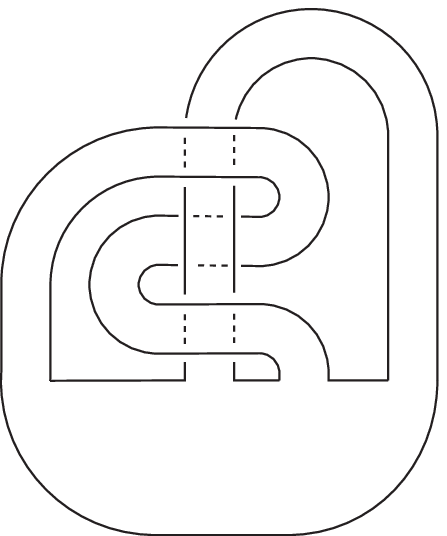}  
\smallbreak
\hskip26mm   F{\tiny IGURE} \ref{Hawaii}.4:
{\bf The resulting nontrivial $n$-knot $K$}

\bigbreak

We can prove that $K$ is nontrivial  
by using Seifert matrices and the Alexander polynomial. 
(See Definition \ref{Missouri} for Seifert matrices, and 
 Definition \ref{square} for the Alexander polynomial.)
We use the fact that $S^p$ and $S^q$ can be `linked' in $S^{p+q+1}$. 
Recall that $p+q+1=n+2$. Note that $S^q$ and $S^p$ are included in $F$ as shown below.

\bigbreak
\hskip4cm\includegraphics[width=6cm]{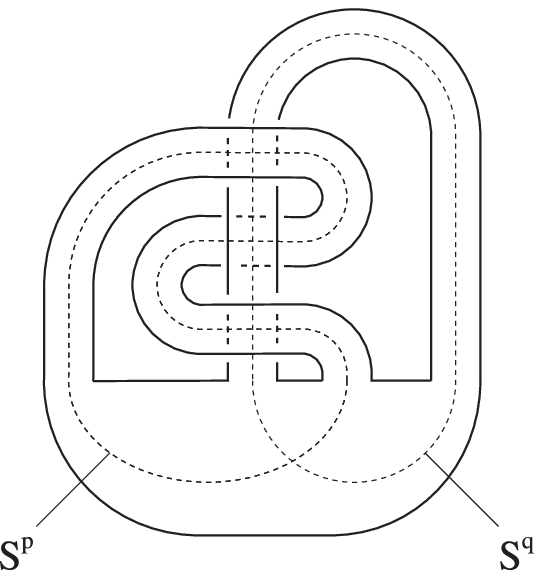}  
\smallbreak
\hskip2cm F{\tiny IGURE} \ref{Hawaii}.5: {\bf $S^p$ and $S^q$ in $F$ whose boundary is $K$} 

\bigbreak

Note that the above operation is done only in an $(n+2)$-ball.  
This operation is an example of the $(p, q)$-pass-moves, 
whose definition we review in Definition \ref{Minesota}. 

Local moves on high dimensional submanifolds are exciting ways of explicit construction of high dimensional figures. They are also generalization of local moves on 1-links. 
They are useful to research link cobordism, knot cobordism, and the intersection of submanifolds (see \cite{Ogasa98n}) etc.
There remain many exciting problems. Some of them are proper in high dimensional case and others are analogous to 1-dimensional one. For example, we do not know a local move on high dimensional knots which is an unknotting operation.  


\bigbreak
Let $p, q\in\N$ and $p+q=n+1$. 
We review the definition of the $(p, q)$-pass-move on $n$-knots,   
which was defined in \cite{Ogasa98n} 
and 
which has been studied in 
\cite{KauffmanOgasa, Ogasa98n,  Ogasa02,  Ogasa04, Ogasa07,   Ogasa09, OgasaT3, OgasaIH}. 
If $p=1$ and $q=1$, 
the $(p,q)$-pass-move on $(p+q-1)$-knots is  the pass-move on 1-knots, 
whose definition we reviewed in \S\ref{Introduction}. 
If $p=1$ and $q=2$, the $(p,q)$-pass-move on $(p+q-1)$-knots is  
the (1,2)-pass-move on 2-knots,  
whose definition we reviewed in \S\ref{Illinois}.

\begin{defn}\label{Minesota}
Let $n,p\in\N$. Let $n+1-p>0.$     
Regard an $(n+2)$-ball $B=D^{n+2}$  
as \newline $D^1\x D^p\x D^{n+1-p}$. 
(See Figure \ref{Hawaii}.7.)  
Let $D^1=[-1,1]=\{t|-1\leqq t\leqq1\}.$ 
Take 
a $p$-ball $D_s^p$ (resp. an $(n+1-p)$-ball $D_s^{n+1-p}$)   
embedded in 
Int$D^p$ (resp.  Int$D^{n+1-p}$.)  
Let 
$S^{p-1}$ (resp. $S^{n-p}$)
denote 
the $(p-1)$-sphere $\partial D_s^p$ (resp. the $(n-p)$-sphere $\partial D_s^{n+1-p}$).  
Let a submanifold 
$\{0\}\x D^p\x D_S^{n+1-p}$  (resp.  $\{0\}\x D_S^p\x D^{n+1-p}$) 
$\subset B$ 
be called $h^p$ (resp. $h^{n+1-p}$). 
We give an orientation to $h^p$ (resp. $h^{n+1-p}$). 
Note that $h^p\cap h^{n+1-p}\neq\phi$. 
Move $h^p$ in $B$ by using an isotopy with keeping $h^p\cap \partial B$, 
 let 
`the resultant submanifold$-B$'   
be put in 
$\{t>0\}\x D^p\x D^{n+1-p}$ (resp.  $\{t<0\}\x D^p\x D^{n+1-p}$), 
and call the submanifold $h_+^p$ (resp. $h_-^p$).
(See Figure \ref{Hawaii}.6.)   
Note that $h_+^p\cap h^{n+1-p}=\phi$, 
that $h_-^p\cap h^{n+1-p}=\phi$, 
and that $h^p_+\cap h^p_-=h^p_+\cap \partial B=h^p_-\cap \partial B$. 
%
%
(Each of 
Figure \ref{Hawaii}.6 
and 
Figure \ref{Hawaii}.7, 
which consists of the two figures (1) and (2),  is a diagram of the $(p, q)$-pass-move, 
where $q=n+1-p$.)

Let $K_+$ and $K_-$ be 
$n$-dimensional closed oriented submanifolds $\subset S^{n+2}$. 
Embed the $(n+2)$-ball $B$ 
in $S^{n+2}$. 
Let $K_+$ and $K_-$ differ only in $B$. 
Let $K_+$ (resp. $K_-$) satisfy the condition  \newline
\hskip39mm$K_+\cap \mathrm{Int}B=
(\partial h^{p}_+-\partial B)\cup(\partial h^{n+1-p}-\partial B)$ \newline  
\hskip28mm 
${\rm(resp.}\hskip1mm K_-\cap \mathrm{Int}B=
(\partial h^{p}_--\partial B)\cup(\partial h^{n+1-p}-\partial B)),$  \newline
where we suppose that there is not 
$h^{p}_-$  
(resp. $h^{p}_+$)   
in $B$.  
Then we say that 
$K_+$ (resp. $K_-$)  is obtained from $K_-$ (resp. $K_+$) 
by one {\it $(p,n+1-p)$-pass-move} in $B$.   
\end{defn}

In Definition \ref{Minesota} we have the following: 
Let $\sharp\in\{+,-\}.$   
there is a Seifert hypersurface $V_\sharp\subset S^{n+2}$ for $K_\sharp$ such that 
$V_\sharp\cap B=h_\sharp^p\cup h^{n+1-p}.$ 
\noindent 
(The idea of the proof is Thom-Pontrjagin construction.)
We say that  
$V_-$ (resp. $V_+$) is obtained from $V_+$ (resp. $V_-$) 
by one {\it $(p,n+1-p)$-pass-move} in $B$.

\bigbreak

In Definition \ref{Minesota}, note the following: 
Let 
$V_0=V_\sharp-\text{Int}B$ \newline
$=\text{`the closure of }(V_\sharp- (h_\sharp^p\cup h^{n+1-p}))\text{ in }{S^{n+2}}'.$  
We can say that 
we attach 
an embedded $(n+1)$-dimensional $p$-handle $h_\#^p$ 
and 
an embedded $(n+1)$-dimensional $(n+1-p)$-handle $h^{n+1-p}$ 
to the submanifold $V_0\subset S^{n+2}$, 
and obtain the submanifold $V_\#\subset S^{n+2}$.

\begin{figure}
\unitlength 0.1in
\begin{picture}(52.00,24.84)(1.50,-32.30)
%
\special{pn 8}%
\special{ar 4488 1824 843 843  0.8934700 6.2831853}%
\special{ar 4488 1824 843 843  0.0000000 0.8617842}%
%
\special{pn 8}%
\special{pa 3829 1312}%
\special{pa 4891 2556}%
\special{fp}%
%
\special{pn 8}%
\special{pa 5155 2326}%
\special{pa 4104 1083}%
\special{fp}%
%
\special{pn 8}%
\special{pa 4891 1083}%
\special{pa 4515 1440}%
\special{fp}%
\special{pa 4168 1806}%
\special{pa 3738 2208}%
\special{fp}%
\special{pa 5147 1320}%
\special{pa 4735 1686}%
\special{fp}%
\special{pa 4369 2071}%
\special{pa 3966 2455}%
\special{fp}%
%
\special{pn 8}%
\special{ar 1652 1858 843 843  5.5899622 6.2831853}%
\special{ar 1652 1858 843 843  0.0000000 5.5598920}%
%
\special{pn 8}%
\special{pa 904 1468}%
\special{pa 1268 1838}%
\special{fp}%
\special{pa 1639 2180}%
\special{pa 2047 2604}%
\special{fp}%
\special{pa 1138 1209}%
\special{pa 1512 1614}%
\special{fp}%
\special{pa 1901 1974}%
\special{pa 2291 2371}%
\special{fp}%
%
\special{pn 8}%
\special{pa 2146 1183}%
\special{pa 917 2254}%
\special{fp}%
%
\special{pn 8}%
\special{pa 1149 2526}%
\special{pa 2377 1445}%
\special{fp}%
%
\special{pn 8}%
\special{pa 884 3002}%
\special{pa 1039 2179}%
\special{dt 0.045}%
\special{sh 1}%
\special{pa 1039 2179}%
\special{pa 1007 2241}%
\special{pa 1029 2231}%
\special{pa 1046 2248}%
\special{pa 1039 2179}%
\special{fp}%
\special{pa 921 2983}%
\special{pa 1258 2426}%
\special{dt 0.045}%
\special{sh 1}%
\special{pa 1258 2426}%
\special{pa 1206 2473}%
\special{pa 1230 2472}%
\special{pa 1241 2493}%
\special{pa 1258 2426}%
\special{fp}%
%
\special{pn 8}%
\special{pa 2402 2801}%
\special{pa 2228 2334}%
\special{dt 0.045}%
\special{sh 1}%
\special{pa 2228 2334}%
\special{pa 2233 2403}%
\special{pa 2247 2384}%
\special{pa 2270 2389}%
\special{pa 2228 2334}%
\special{fp}%
\special{pa 2393 2783}%
\special{pa 1853 2361}%
\special{dt 0.045}%
\special{sh 1}%
\special{pa 1853 2361}%
\special{pa 1893 2418}%
\special{pa 1895 2394}%
\special{pa 1918 2386}%
\special{pa 1853 2361}%
\special{fp}%
\put(22.1900,-30.6600){\makebox(0,0)[lb]{$S^{p-1}\x D^{n+1-p}$}}%
\put(3.0700,-31.6600){\makebox(0,0)[lb]{$D^p\x S^{n-p}$}}%
\put(21.8200,-33.0500){\makebox(0,0)[lb]{$=\overline{\partial h^{n+1-p}-\partial B}$}}%
\put(11.9400,-9.5300){\makebox(0,0)[lb]{$B\cap K_+$}}%
\put(40.6700,-9.1600){\makebox(0,0)[lb]{$B\cap K_-$}}%
\put(1.5000,-34.0000){\makebox(0,0)[lb]{$=\overline{\partial h^p_+-\partial B}$}}%
%
\special{pn 8}%
\special{pa 1100 1450}%
\special{pa 1111 1419}%
\special{pa 1121 1389}%
\special{pa 1128 1358}%
\special{pa 1130 1326}%
\special{pa 1128 1295}%
\special{pa 1122 1263}%
\special{pa 1112 1232}%
\special{pa 1100 1201}%
\special{pa 1087 1171}%
\special{pa 1071 1143}%
\special{pa 1053 1116}%
\special{pa 1033 1092}%
\special{pa 1009 1072}%
\special{pa 982 1055}%
\special{pa 954 1040}%
\special{pa 924 1026}%
\special{pa 893 1015}%
\special{pa 862 1004}%
\special{pa 831 1000}%
\special{pa 800 1007}%
\special{pa 770 1020}%
\special{sp}%
\put(4.1000,-11.9000){\makebox(0,0)[lb]{$h^{n+1-p}$}}%
\put(32.6000,-11.8000){\makebox(0,0)[lb]{$h^{n+1-p}$}}%
%
\special{pn 8}%
\special{pa 4040 1360}%
\special{pa 4035 1327}%
\special{pa 4028 1296}%
\special{pa 4018 1266}%
\special{pa 4002 1238}%
\special{pa 3982 1214}%
\special{pa 3961 1189}%
\special{pa 3941 1165}%
\special{pa 3921 1140}%
\special{pa 3902 1115}%
\special{pa 3883 1089}%
\special{pa 3864 1062}%
\special{pa 3844 1035}%
\special{pa 3821 1012}%
\special{pa 3794 998}%
\special{pa 3761 994}%
\special{pa 3729 1000}%
\special{pa 3700 1015}%
\special{pa 3670 1020}%
\special{sp}%
%
\special{pn 8}%
\special{pa 2190 1380}%
\special{pa 2199 1349}%
\special{pa 2209 1318}%
\special{pa 2219 1288}%
\special{pa 2231 1259}%
\special{pa 2245 1230}%
\special{pa 2261 1202}%
\special{pa 2280 1176}%
\special{pa 2300 1151}%
\special{pa 2322 1127}%
\special{pa 2346 1106}%
\special{pa 2373 1089}%
\special{pa 2404 1078}%
\special{pa 2435 1069}%
\special{pa 2464 1058}%
\special{pa 2488 1040}%
\special{pa 2508 1016}%
\special{pa 2524 987}%
\special{pa 2539 956}%
\special{pa 2550 930}%
\special{sp}%
\put(25.0000,-9.5000){\makebox(0,0)[lb]{$h^p_+$}}%
\put(53.5000,-12.4000){\makebox(0,0)[lb]{$h^p_-$}}%
%
\special{pn 8}%
\special{pa 5000 1290}%
\special{pa 5028 1274}%
\special{pa 5054 1255}%
\special{pa 5078 1234}%
\special{pa 5098 1208}%
\special{pa 5117 1180}%
\special{pa 5140 1160}%
\special{pa 5170 1152}%
\special{pa 5203 1145}%
\special{pa 5234 1140}%
\special{pa 5265 1145}%
\special{pa 5297 1150}%
\special{pa 5320 1150}%
\special{sp}%
\end{picture}%

\vskip5mm
\hskip5mm\text{ F{\tiny IGURE} \ref{Hawaii}.6: 
{\bf 
The $(p, n+1-p)$-pass-move on an $n$-dimensional closed submanifold  \newline
 }}

\hskip-30mm\text{ {\bf 
$\subset S^{n+2}$. 
Note $B=B^{n+2}=D^{n+2}\subset S^{n+2}$.   
}}

\end{figure}

\begin{figure}
\begin{center}
\includegraphics[width=107mm]{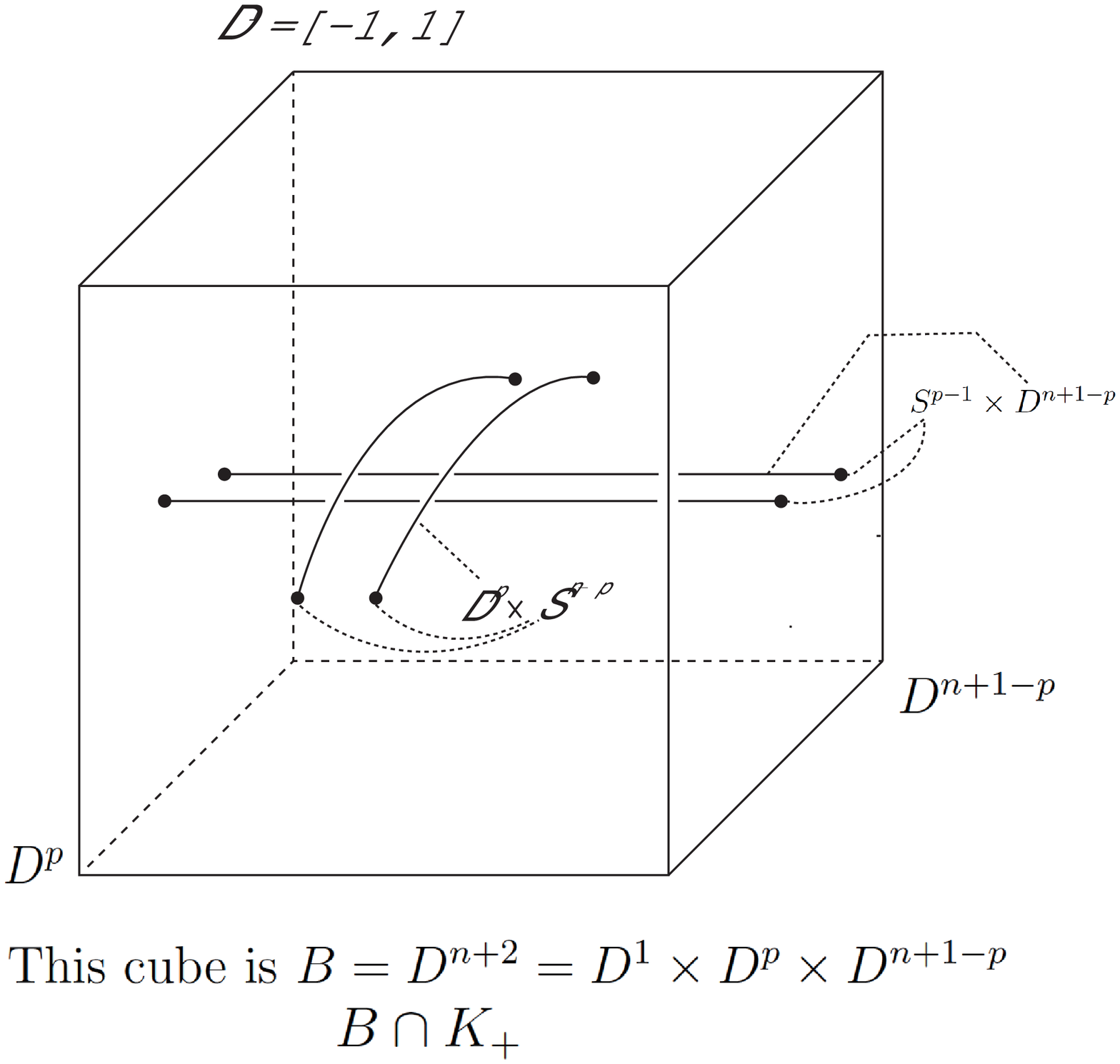}

\text{
F{\tiny IGURE} \ref{Hawaii}.7.(1): {\bf The $(p,n+1-p)$-pass-move}
}

\vskip10mm
\includegraphics[width=70mm]{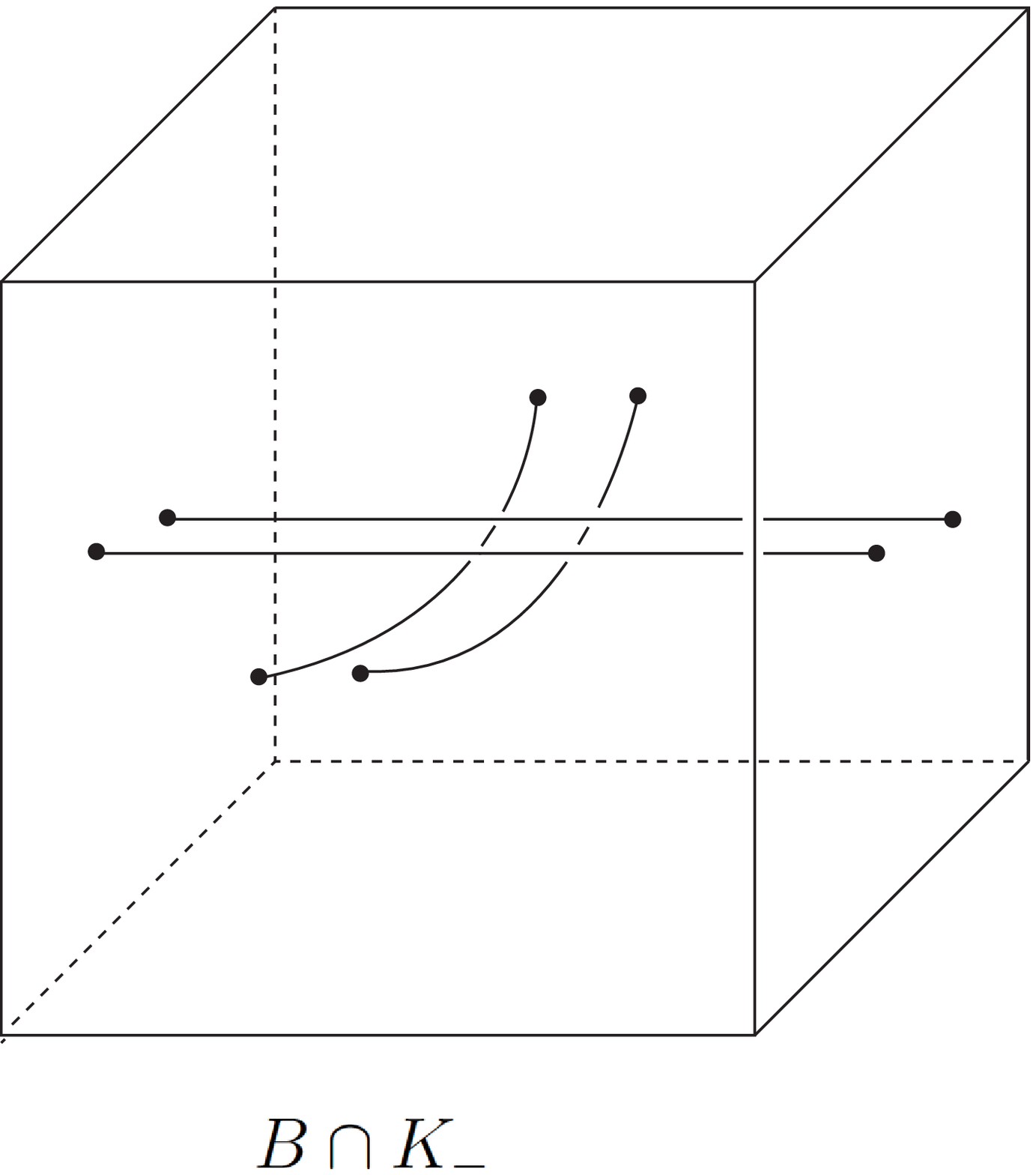}

\vskip2mm
\text{F{\tiny IGURE} \ref{Hawaii}.7.(2): {\bf The $(p,n+1-p)$-pass-move} }
\end{center}
\end{figure}

\begin{defn}\label{Montana}
Let $p,q,p+q-1\in\N.$
Let $K$ be a $(p+q-1)$-knot$\subset S^{p+q+1}$ 
which is $(p,q)$-pass-move-equivalent to the trivial $(p+q-1)$-knot. 
The {\it $(p,q)$-pass-move-unknotting-number} of $K$ is 
the minimal number of $(p,q)$-pass-moves which we change $K$ to the trivial $(p+q-1)$-knot by. 
\end{defn}

\begin{prop}\label{Nebraska}
Let $p,q,p+q-1\in\N.$
There is a $(p+q-1)$-knot 
whose $(p,q)$-pass-move-unknotting-number is one. 
\end{prop}

\noindent
{\bf Proof of Proposition \ref{Nebraska}.}  
See the nontrivial $(p+q-1)$-knot  $K$ in Figures \ref{Hawaii}.4 of this section: 
A Seifert hypersurface $V$ for $K$ is diffeomorphic to $S^p\x S^q-{\rm Int}B^{p+q}$. 
We supposed the following: 
$x$ (resp. $y$) is a generator of $H_p(V;\Z)$ (resp. $H_p(V;\Z)$). 
The intersection matrix associated with the base $\{x\}$ and  $\{y\}$
is a $1\x1$-matrix $(1)$.  
The Seifert matrix associated with the base $\{x\}$ and  $\{y\}$
is a $1\x1$-matrix $(2)$.  

Hence the $(p,q)$-pass-move-unknotting-number of $K$ is$\geqq1$.

$K$ is obtained from the trivial $(p+q-1)$-knot by 
one $(p,q)$-pass-move as drawn in Figures  \ref{Hawaii}.2-4. 

Therefore the $(p,q)$-pass-move-unknotting-number of $K$ is one. 
\qed\bigbreak

We consider the following problem. 

\begin{prob}\label{Louisiana}   
Let $k\in\N\cup\{0\}.$
\smallbreak\noindent
$(1)$ Is there a 
$(2k+1, 2k+2)$-pass-move-unknotting-number-two $(4k+2)$-knot? 

\smallbreak\noindent
$(2)$ For any natural number $n$, 
is there a 
$(4k+2)$-knot 
whose $(2k+1, 2k+2)$-pass-move-unknotting-number is$>n$?  
\end{prob}

We give a positive answer to Problem \ref{Louisiana}.(1) (resp. \ref{Louisiana}.(2)). 
The answers make one of our main theorems.   

\begin{thm}\label{Indiana}  
Let $k\in\N\cup\{0\}.$
\smallbreak\noindent
$(1)$ There is a $(2k+1, 2k+2)$-pass-move-unknotting-number-two $(4k+2)$-knot. 

\smallbreak\noindent
$(2)$ For any natural number $n$, there is a $(4k+2)$-knot 
whose $(2k+1, 2k+2)$-pass-move-unknotting-number is$>n$. 
\end{thm}

We  consider the following problem. 

\begin{prob}\label{Maine}   
Let $k\in\N\cup\{0\}.$
\smallbreak\noindent
$(1)$ Is there a 
$(2k+1, 2k+1)$-pass-move-unknotting-number-two $(4k+1)$-knot? 

\smallbreak\noindent
$(2)$ For any natural number $n$, 
is there a 
$(4k+1)$-knot 
whose $(2k+1, 2k+1)$-pass-move-unknotting-number is$>n$?  
\end{prob}

We give a positive answer to Problem \ref{Maine}.(1) (resp. \ref{Maine}.(2)). 
The answers make one of our main theorems.   

\begin{thm}\label{Iowa}  
Let $k\in\N\cup\{0\}.$
\smallbreak\noindent
$(1)$ There is a $(2k+1, 2k+1)$-pass-move-unknotting-number-two $(4k+1)$-knot. 

\smallbreak\noindent
$(2)$ For any natural number $n$, there is a $(4k+1)$-knot 
whose $(2k+1, 2k+1)$-pass-move-unknotting-number is$>n$. 
\end{thm}

\bigbreak
\section{Proof of Theorem \ref{Indiana}}\label{Maryland}






-


\begin{prop}\label{hMississippi}
Let $k\in\N\cup\{0\}.$
Let $K$  be a $(4k+2)$-knot $\subset S^{4k+4}$   
whose $(2k+1,2k+2)$-pass-move-unknotting-number is one.  
Let $M_3(K)$ be the 3-fold branched covering space of $S^{4k+4}$ along $K$. 
Then there are three elements$\in H_{2k+1}(M_3(K);\Z)$ 
which generate $H_{2k+1}(M_3(K);\Z)$. 
\end{prop}

\begin{figure}
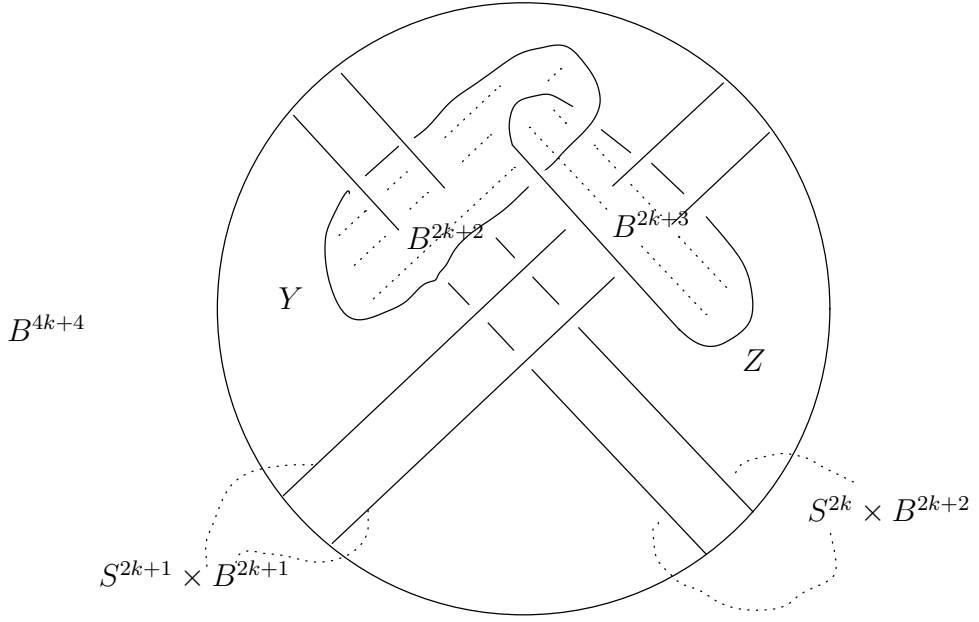
 
\input hsurgery.tex 

\smallbreak
\caption{{\bf The (2k+1,2k+2)-pass-move carried out by surgeries}\label{neko}} 
\smallbreak
\end{figure}

\noindent{\bf Proof of Proposition \ref{hMississippi}.}
Take a $(4k+4)$-ball $B^{4k+4}\subset S^{4k+4}$ where 
we carry out the  $(2k+1,2k+2)$-pass-move which changes $K$ into $T$. 
See Figure \ref{neko}. 
Note that 
$K\cap B^{4k+4}\newline 
=(S^{2k}\x B^{2k+2})\amalg(S^{2k+1}\x B^{2k+1})$.  
Take a $(2k+2)$-ball $B^{2k+2}$ in the $(4k+4)$-ball $B^{4k+4}$ 
such that 
$B^{2k+2}\cap(S^{2k}\x B^{2k+2})$ is the $2k$-sphere trivially embedded in $B^{2k+2}$. 
and such that 
$B^{2k+2}\cap(S^{2k+1}\x B^{2k+1})=\phi$. 
Call $\partial B^{2k+2}$, $Y$. 
Take a $(2k+3)$-ball $B^{2k+3}$ in the $(4k+4)$-ball $B^{4k+4}$ 
such that 
$B^{2k+3}\cap(S^{2k+1}\x B^{2k+1})$ is the $(2k+1)$-sphere trivially embedded in $B^{2k+3}$ 
and such that 
$B^{2k+3}\cap(S^{2k}\x B^{2k+2})=\phi$.   
Call $\partial B^{2k+3}$, $Z$.  
Suppose that the linking number of $Y$ and $Z$ is one. 
Attach a $(4k+5)$-dimensional $(2k+2)$-(resp. $(2k+3)$-)handle to $B^{4k+4}$ 
along $Y$ (resp. $Z$) with the trivial framing.
Note that these two handles are attached to $S^{4k+4}$ on time. 
Carry out surgeries by using these two handles on $S^{4k+4}$.   
Then the new manifold which we obtain is the $(4k+4)$-sphere again, and call it $S^{'4k+4}$.  
Furthermore the new submanifold $\subset S^{'4k+4}$ which is made from $K$  
is the trivial $(4k+2)$-knot $T$. 
 
Note that now we have a compact oriented $(4k+5)$-dimensional manifold $W$ with 
a handle decomposition 

\smallbreak\hskip3mm
$W=(S^{4k+4}\x[0,1])\cup$(a $(4k+5)$-dimensional $(2k+2)$-handle)

\hskip14mm $\cup$(a $(4k+5)$-dimensional $(2k+3)$-handle)$\cup(S^{'4k+4}\x[0,1])$. 
\smallbreak

Note that $\partial W=(S^{4k+4}\x\{0\})\amalg(S^{'4k+4}\x\{1\})$.   
Note that there is an embedding map $f:S^{4k+2}\x [0,1]\hookrightarrow W$ with the following properties: 

\smallbreak
\noindent
(1) 
$f(S^{4k+2}\x [0,1])\cap(S^{4k+4}\x\{0\})$ 
is  $f(S^{4k+2}\x\{0\})$. 
$f(S^{4k+2}\x [0,1])\cap(S^{'4k+4}\x\{1\})$ 
is  $f(S^{4k+2}\x\{1\})$. 

$f$ is transverse to $\partial W.$

\smallbreak
\noindent
(2) 
 $f(S^{4k+2}\x\{0\})$ in $(S^{4k+4}\x\{0\})$ is $K$. 

\hskip2mm$f(S^{4k+2}\x\{1\})$ in $(S^{'4k+4}\x\{1\})$ is $T$. 

\bigbreak

Take a 3-fold branched covering space $\widetilde{W}$ of $W$ along $f(S^{4k+2}\x [0,1])$. 
Note the $(2k+1)$-sphere which is the core of the attaching part of the $(2k+2)$-handle in 
the above handle decomposition of $W$. 
The $(2k+1)$-sphere is null-homologous in $S^{4k+4}-N(K)$, 
where $N(K)$ is the tubular neighborhood of $K$ in $S^{4k+4}$.    
Therefore we obtain a compact oriented $(4k+5)$-dimensional manifold $\widetilde{W}$ with 
a handle decomposition 

\smallbreak\hskip3mm
$\widetilde{W}=
(M_3(K)\x[0,1])$

\quad\quad$\cup$
(three $(4k+5)$-dimensional $(2k+2)$-handles, $h^{2k+2}_1$, $h^{2k+2}_2$, and $h^{2k+2}_3$)

\quad\quad$\cup$
(three $(4k+5)$-dimensional $(2k+3)$-handle, $h^{2k+3}_1$, $h^{2k+3}_2$, and $h^{2k+3}_3$)

\quad\quad$\cup(S^{'4k+4}\x[0,1])$. 
\smallbreak
 
\noindent
Here, note that the 3-fold branched covering space of $S^{'4k+4}$ along $T$ 
is the standard $(4k+4)$-sphere, 
and call it $S^{'4k+4}$ again.

We prove that $H_{2k+1}(\widetilde{W};\Z)\cong0$. 
{\it Reason.}
Take the dual handle decomposition 

\smallbreak
\hskip3mm
$\widetilde{W}=(S^{'4k+4}\x[0,1])$

\quad\quad$\cup$(three $(4k+5)$-dimensional $(2k+2)$-handles, 
$\overline{h^{2k+2}_1}$, 
$\overline{h^{2k+2}_2}$, 
$\overline{h^{2k+2}_3}$)

\quad\quad$\cup$(three $(4k+5)$-dimensional $(2k+3)$-handles, 
$\overline{h^{2k+3}_1}$, 
$\overline{h^{2k+3}_2}$, 
$\overline{h^{2k+3}_3}$)

\quad\quad$\cup(M_3(K)\x[0,1])$, 
\smallbreak

\noindent
of the above handle decomposition, where $\overline{h^*_\#}$     
is the dual handle of $h^{4k+5-*}_\#$. 
Take a manifold $Q_S$ which is represented by the sub-handle-decomposition  

$Q_S$

$=(S^{'4k+4}\x[0,1])\cup$(the three $(4k+5)$-dimensional $(2k+2)$-handles, 
$\overline{h^{2k+2}_1}$, 
$\overline{h^{2k+2}_2}$, 
$\overline{h^{2k+2}_3}$)

\noindent
of the dual handle decomposition of $\widetilde{W}$. 
Since $H_{2k+1}(S^{'4}\x[0,1];\Z)\cong0$,  we have 
 $H_{2k+1}(Q_S;\Z)\cong0$. 
Recall that if we attach $(2k+3)$-handles to a manifold $E$ and we obtain a new manifold $E'$, 
then $H_{2k+1}(E;\Z)\cong H_{2k+1}(E';\Z)$. 


Therefore the manifold $R_S$ 
which is represented by the sub-handle-decomposition 

$R_S=$

$(M_3(K)\x[0,1])\cup($three $(4k+5)$-dimensional $(2k+2)$-handles, $h^{2k+2}_1$, $h^{2k+2}_2$, and $h^{2k+2}_3$)

\noindent 
of the above handle decomposition satisfies 
the condition $H_{2k+1}(R_S;\Z)\cong0.$

Therefore 
 the cores of the attaching parts of $h^{2k+2}_1$, $h^{2k+2}_2$ and $h^{2k+2}_3$ 
generate \newline$H_{2k+1}(M_3(K);\Z)$.     

This completes the proof of Proposition \ref{hMississippi}. \qed\bigbreak

In a similar fashion we can prove the following. 
\begin{prop}\label{hgMississippi}
Let $k\in\N\cup\{0\}.$ Let $n\in\N.$
Let $K$  be a $(4k+2)$-knot $\subset S^{4k+4}$   
whose $(2k+1,2k+2)$-pass-move-unknotting-number of $K$ is$\leqq n$.  
Let $M_3(K)$ be the 3-fold branched covering space of $S^{4k+4}$ along $K$. 
Then there are $3n$ elements$\in H_{2k+1}(M_3(K);\Z)$ 
which generate $H_{2k+1}(M_3(K);\Z)$. 
\end{prop}

\begin{cla}\label{Nevada}  
Let $k\in\N\cup\{0\}.$  
There is a $(4k+2)$-knot $P\subset S^{4k+4}$ as follows. 

\smallbreak\noindent 
$(1)$    
A Seifert hypersurface $V$ for $P$ is diffeomorphic to \newline 
$((S^{2k+1}\x S^{2k+2})\sharp(S^{2k+1}\x S^{2k+2}))-\text{open} B^{4k+3}$.

For  an ordered set  $(x_1,x_2)$ of basis of $H_{2k+1}(V;\Z)\cong\Z\oplus\Z$ 
and  an ordered set  $(y_1,y_2)$ of basis of $H_{2k+2}(V;\Z)\cong\Z\oplus\Z$,  
the intersection matrix $(x_k\cdot y_l)$ $(k, l\in\{1,2\})$   
on $H_{2k+1}(V;\Z)$ $($resp. $H_{2k+2}(V;\Z))$ is 
$
\begin{pmatrix}
0&1\\
-1&0
\end{pmatrix}. 
$
We can suppose that Poincar\'e dual of 
$x_1$ $($resp. $x_2)$ is $y_2$ $($resp. $-y_1).$

\smallbreak\noindent 
$(2)$   
The Seifert matrix $({\rm lk}(x_k, y_l^+))$   for $P$  
is  
$X=
\begin{pmatrix}
-1&1\\
0&-1
\end{pmatrix}. 
$
\end{cla}

\noindent
{\bf Note.} 
The negative Seifert matrix related to $X$ is the transposed matrix of $X$. 
Recall that 
a $(2k+1)$-positive Seifert matrix of $(4k+2)$-knot 
is not the transposed matrix of its related negative Seifert matrix in general 
(see Definitions \ref{Missouri} and \ref{square}, and Propositions \ref{tsuika} and \ref{Mars}). 
Note that we have the following: 

\smallbreak   
$
X^{-1}=
\begin{pmatrix}
-1&-1\\
0&-1
\end{pmatrix}, 
$
$
^t\hskip-1mm X X^{-1}=
\begin{pmatrix}
1&1\\
-1&0
\end{pmatrix}, 
$
$
(^t\hskip-1mm X X^{-1})^2=
\begin{pmatrix}
0&1\\
-1&-1
\end{pmatrix}, 
$
$
(^t\hskip-1mm X X^{-1})^3=
\begin{pmatrix}
-1&0\\
0&-1
\end{pmatrix}.  
$

\begin{figure}
\includegraphics[width=10cm]{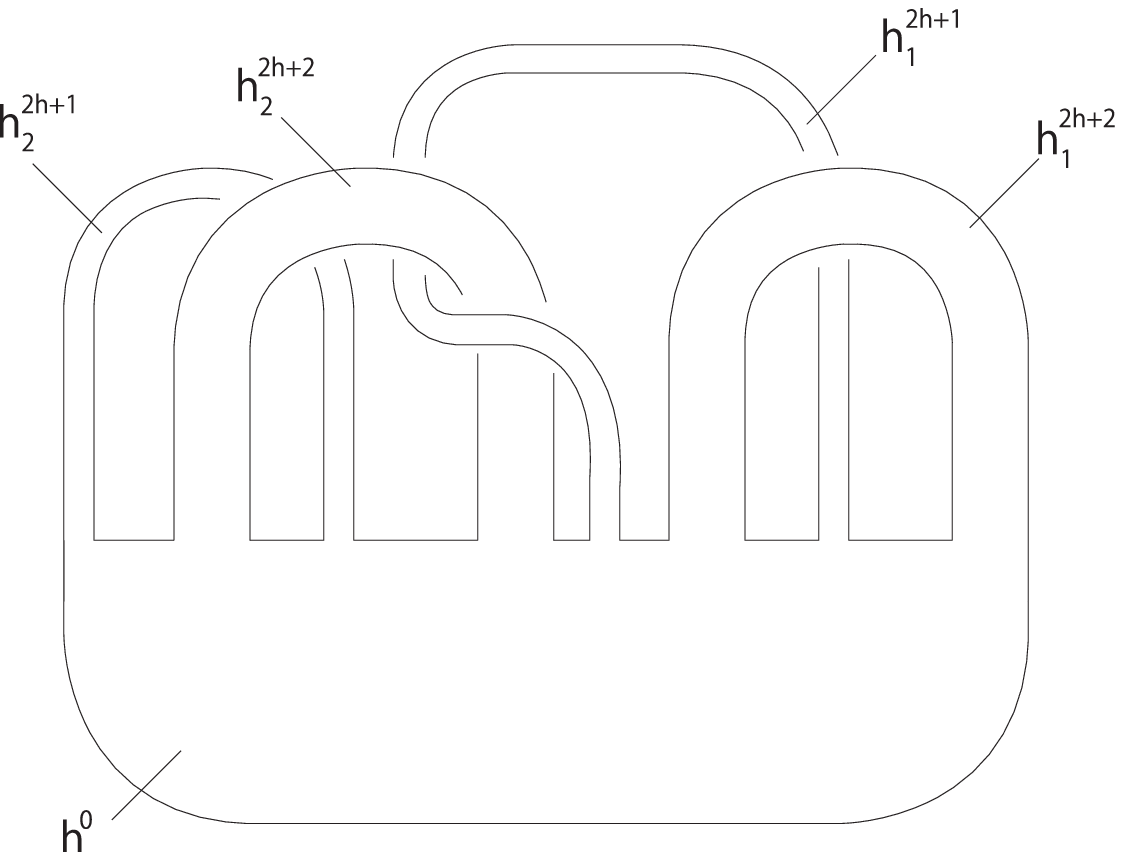}

\smallbreak
$h^0$ denotes a $(4k+2)$-dimensional 0-handle. 

$h^{*}_i (i=1,2 \rm{and} *=2k+1, 2k+2)$ denotes a $(4k+2)$-dimensional $*$-handle. 

$h^{2k+1}_i (i=1,2)$ corresponds $x_i$.
 
$h^{2k+2}_i (i=1,2)$ corresponds $y_i$.
 
\caption{{\bf A Seifert hypersurface $V_T$ for a $(4k+2)$-knot $T$}\label{shishi} }
\end{figure}

\bigbreak
\noindent{\bf Proof of Claim \ref{Nevada}.}
See Figure \ref{shishi}. 
Embed $((B^{2k+2}\x S^{2k+2})\natural(B^{2k+2}\x S^{2k+2}))$ in $S^{4k+4}$. 
Note that its boundary is diffeomorphic to 
$((S^{2k+1}\x S^{2k+2})\sharp(S^{2k+1}\x S^{2k+2}))$.  
Remove an open  $B^{4k+3}$ from it. 
We can suppose that this \newline 
$((S^{2k+1}\x S^{2k+2})\sharp(S^{2k+1}\x S^{2k+2}))-\text{open} B^{4k+3}$
is a Seifert hypersurface $V_T$ for the trivial $(4k+2)$-knot. 
We can take  
 an ordered set of basis $(x_1,x_2)$ (resp. $(y_1, y_2)$) of $H_{2k+1}(V_T;\Z)$ (resp. $H_{2k+2}(V_T;\Z)$) 
%
which satisfies (1) of Claim \ref{Nevada}.  
Furthermore we can suppose that 
the Seifert matrix $({\rm lk}(x_k, y_l^+))$  associated with $V_T$ and this ordered set of basis 
is  
$X=
\begin{pmatrix}
0&1\\
0&0
\end{pmatrix}. 
$

Note that 
$V_T$ has a handle decomposition 

\smallbreak
(0-handle$h^0$)$\cup$(two $(4k+3)$-dimensional $(2k+1)$-handles, $h^{2k+1}_1$ and $h^{2k+1}_2$)

$\cup$(two $(4k+3)$-dimensional $(2k+2)$-handles, $h^{2k+2}_1$ and $h^{2k+2}_2$). 
\smallbreak

\noindent
By two times of $(2k+1,2k+2)$-pass-move     
we change the submanifold $V_T$ into the submanifold $V$ 
so that $V$ satisfies (1) and (2) of Claim \ref{Nevada}.

This completes the proof of Claim \ref{Nevada}.    
\qed\bigbreak

\noindent{\bf Note.}  
We can say that   
the $(4k+2)$-knot $P$ in Claim \ref{Nevada} 
is the knot product of 
the 2-knot $P$ in \S\ref{Kansas} and 
$k$ copies of the Hopf link. 
See \cite{KauffmanOgasa} for the knot product.

\bigbreak  
In Proof of Claim \ref{Nevada} we also prove that 
the $(2k+1, 2k+2)$-pass-move-unknotting-number of $P$ in Proof of Claim \ref{Nevada}  
is$\leqq2$.  
(Another proof is given by using Main Theorem 2.6 of \cite{KauffmanOgasaII}.)

Let $M_3(P)$ be the 3-fold branched covering space of $S^{4k+4}$ along $P$. 

By Proposition \ref{Alabama} 
we have $H_{2k+1}(M_3(P);\Z)\cong\Z_2\oplus\Z_2\oplus\Z_2\oplus\Z_2$. 
Hence we need no less than four generators in order to generate 
$H_{2k+1}(M_3(P);\Z)$. 

Suppose that 
the $(2k+1, 2k+2)$-pass-move-unknotting-number of $P$ is$\leqq1$.
By Proposition \ref{hMississippi},  
we can take three generators in order to generate $H_{2k+1}(M_3(P);\Z)$. 
We arrived at a contradiction.

Therefore the $(2k+1, 2k+2)$-pass-move-unknotting-number of $P$ is$\geqq2$.

Therefore the $(2k+1, 2k+2)$-pass-move-unknotting-number of $P$ is two.

This completes the proof of Theorem \ref{Indiana}.(1). 
\bigbreak

Let $n\in\N$. 
Let $m\in\N$ and $\frac{2m}{3}>n$.  
Let $\#^mP$ be the connected-sum of $m$-copies of $P$.

Since  $P$ is $(2k+1, 2k+2)$-pass-move equivalent to the trivial $(4k+2)$-knot, 
$\#^mP$ is $(2k+1, 2k+2)$-pass-move equivalent to the trivial $(4k+2)$-knot.

Let $N_3(\#^mP)$ be the 3-fold branched covering space of $S^4$ along $\#^mP$. 
By Proposition \ref{Alabama} 
we have $H_1(M_3(\#^mP);\Z)\cong\oplus^{2m}\Z_2$. 
Hence we need no less than $2m$ generators in order to generate 
$H_1(M_3(P);\Z)$.

Suppose that 
the $(2k+1, 2k+2)$-pass-move-unknotting-number of $\#^mP$ is$\leqq n$.
By Proposition \ref{hgMississippi} 
$H_1(M_3(\#^mP);\Z)$ can take  $3n$ generators. 
Since $2m>3n$, we arrived at a contradiction. 

Therefore the $(2k+1, 2k+2)$-pass-move-unknotting-number of $\#^mP$ is$>n$.  

This completes the proof of Theorem \ref{Indiana}.(2). 

This completes the proof of Theorem \ref{Indiana}.\qed

\bigbreak
\section{Proof of Theorem \ref{Iowa}}\label{Massachusetts}
\begin{prop}\label{hsix}  
Let $k\in\N\cup\{0\}.$  
Let $J$ be a $(4k+1)$-knot $\subset S^{4k+3}$ 
whose \newline $(2k+1, 2k+1)$-pass-move-unknotting-number is one. 
Let $N_3(J)$ be a 3-fold branched covering space of $S^{4k+3}$ along $J$. 
Then there are six elements$\in H_{2k+1}(N_3(J);\Z)$ 
which generate $H_{2k+1}(N_3(J);\Z)$. 
\end{prop}

\begin{figure}
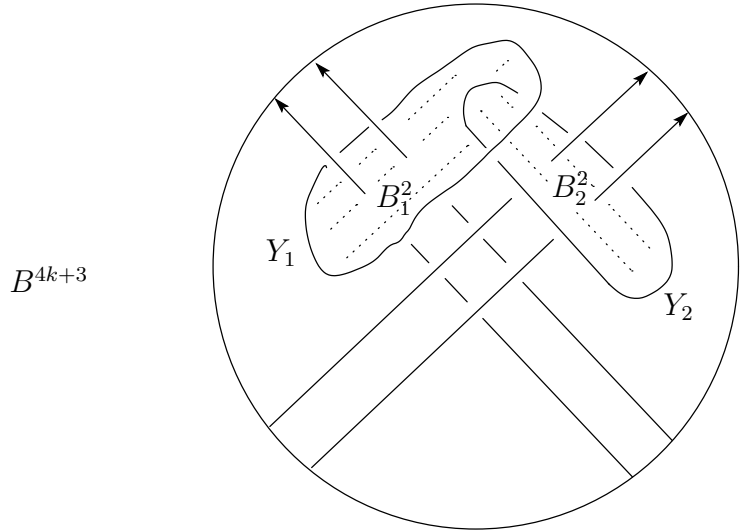
 
\vskip-11mm
\input pass.tex

\hskip3cm \caption{{\bf The $(2k+1,2k+1)$-pass-move carried out by surgeries}\label{kame}}
\smallbreak
\end{figure}

\noindent {\bf Proof of Proposition \ref{hsix}.}   
Take a $(4k+3)$-ball $B^{4k+3}\subset S^{4k+3}$ 
where we carry out the pass-move which changes $J$ into $T$. 
See Figure \ref{kame}.  
Note that $J\cap B^{4k+3}$ is regarded as 
$(S^{2k}\x B^{2k+1})\amalg(S^{2k}\x B^{2k+1})$. 
Call one of the two $S^{2k}\x B^{2k+1}$, $A_1$, and 
the other $A_2$.    
Take two $(2k+2)$-balls, $B^{2k+2}_1$ and $B^{2k+2}_2$ 
in the $(4k+3)$-ball $B^{4k+3}$ such that 
$B^{2k+2}_i\cap A_i$ is a $2k$-sphere trivially embedded in  $B^{2k+2}_i$ 
and such that $B^{2k+2}_i\cap A_j=\phi$  
($i=1,2$, and $i\neq j$).    
Call $\partial B^{2k+2}_i$, $Y_i$ ($i=1,2$).   
Suppose that the linking number of $Y_1$ and $Y_2$ is one. 
Attach a $(4k+4)$-dimensional $(2k+2)$-handle to $B^3$ 
along $Y_i$ with the trivial framing ($i=1,2$).  
Note that these two handles are attached to $S^{4k+3}$ on time. 
Carry out surgeries by using these two handles on $S^{4k+3}$.   
Then the new manifold which we obtain is the $(4k+3)$-sphere again, and call it $S^{'4k+3}$.  
Furthermore the new submanifold$\subset S^{'4k+3}$ which is made from $J$ 
is the trivial $(4k+1)$-knot $T$. 
 
Note that we now have a compact oriented $(4k+4)$-dimensional manifold $U$ with 
a handle decomposition 

\smallbreak\hskip3mm
$U=(S^{4k+3}\x[0,1])\cup$(two $(4k+3)$-dimensional $(2k+2)$-handles)
$\cup(S^{'4k+3}\x[0,1])$. 
\smallbreak

Note that $\partial U=(S^{4k+3}\x\{0\})\amalg(S^{'4k+3}\x\{1\})$. 
There is an embedding map \newline
$f:S^{4k+1}\x [0,1]\hookrightarrow U$ with the following properties: 

\smallbreak
\noindent
(1) 
$f(S^{4k+1}\x [0,1])\cap(S^{4k+3}\x\{0\})$ 
is  $f(S^{4k+1}\x\{0\})$. 
$f(S^{4k+1}\x [0,1])\cap(S^{'4k+3}\x\{1\})$ 
is  $f(S^{4k+1}\x\{1\})$. 

$f$ is transverse to $\partial U.$

\smallbreak
\noindent
(2) 
$f(S^{4k+1}\x\{0\})$ in $S^{4k+3}\x\{0\}$ is $J$. 

$f(S^{4k+1}\x\{1\})$ in $S^{'4k+3}\x\{1\}$ is $T$. 

\bigbreak

Take a 3-fold branched covering space $\widetilde{U}$ 
of $U$ along $f(S^{4k+1}\x [0,1])$. 
Note the \newline
$(2k+1)$-sphere 
which is the core of the attaching part of each of the two $(2k+2)$-handles in 
the above handle decomposition of $U$. 
Each of the two $(2k+1)$-spheres is null-homologous in $S^{4k+3}-N(J)$, 
where $N(J)$ is the tubular neighborhood of $J$ in $S^{4k+3}$.    
Therefore we obtain a compact oriented $(4k+4)$-dimensional manifold $\widetilde{U}$ with 
a handle decomposition 

\smallbreak\hskip3mm
$\widetilde{U}=(N_3(J)\x[0,1])\cup$
(six $(4k+4)$-dimensional $(2k+2)$-handles, $h^{2k+2}_1$,...,$h^{2k+2}_6$)

\quad\quad$\cup(S^{'4k+3}\x[0,1])$.  
\smallbreak
 
\noindent 
Here, note that a 3-fold branched covering space of $S^{4k+3}$ along $T$ is 
the standard \newline $(4k+3)$-sphere, 
and call it $S^{'4k+3}$ again.

We prove that 
$H_{2k+1}(\widetilde{U};\Z)\cong0$. 
{\it Reason.}
Take the dual handle decomposition 

\smallbreak\hskip3mm
$\widetilde{U}=(S^{'4k+3}\x[0,1])\cup$
(six $(4k+4)$-dimensional $(2k+2)$-handles, $\overline{h^{2k+2}_1}$,,...,$\overline{h^{2k+2}_6}$)

\quad\quad
$\cup(N_3(J)\x[0,1])$, 
\smallbreak

\noindent
of the above handle decomposition, where $\overline{h^{2k+2}_\#}$ 
is the dual handle of $h^{2k+2}_\#$. 
Since \newline 
$H_{2k+1}(S^{'3}\x[0,1] ;\Z)\cong0$,    
we have 
$H_{2k+1}(\widetilde{U};\Z)\cong0$.


Therefore  the cores of the attaching parts of 
$h^{2k+2}_1$,...,$h^{2k+2}_6$ generate  
   $H_{2k+1}(N_3(J);\Z)$.   

This completes the proof of Proposition \ref{hsix}.  \qed  \bigbreak 

In a similar way, we can prove the following. 

\begin{prop}\label{hgsix}
Let $k\in\N\cup\{0\}.$  
Let $n\in\N$. 
Let $J$ be a $(4k+1)$-knot $\subset S^{4k+3}$ 
whose $(2k+1, 2k+1)$-pass-move-unknotting-number is$\leqq n$.
Let $N_3(J)$ be a 3-fold branched covering space of $S^{4k+3}$ along $J$. 
Then there are $6n$ elements$\in H_{2k+1}(N_3(J);\Z)$ 
which generate $H_{2k+1}(N_3(J);\Z)$. 
\end{prop}

\begin{figure}
\includegraphics[width=7cm]{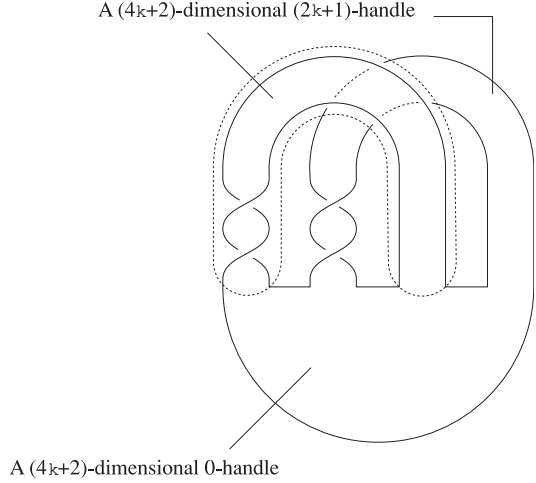}

\caption{{\bf A $(4k+1)$-knot $R$}\label{kujira}}
\end{figure}

Let $R$ be a $(4k+1)$-knot$\subset S^{4k+3}$ 
whose Seifert hypersurface $V$ is diffeomorphic to \newline 
$(S^{2k+1}\x S^{2k+1})-B^{4k+2}$. 
See Figure \ref{kujira}. 
Take  an ordered set  $(x_1,x_2)$ of basis of $H_{2k+1}(V;\Z)\cong\Z\oplus\Z$ 
with the following properties: 

\smallbreak\noindent (1)  
The intersection matrix $(x_k\cdot x_l)$ ($k, l\in\{1,2\}$)  
on $H_{2k+1}(V;\Z)$ is 
$
\begin{pmatrix}
0&1\\
-1&0
\end{pmatrix}. 
$
Note that Poincar\'e dual of $x_1$ is $y_2$.

\smallbreak\noindent (2)  
The Seifert matrix $({\rm lk}(x_k, y_l^+))$   for $R$ is    
$X=
\begin{pmatrix}
-1&1\\
0&-1
\end{pmatrix}
$

\smallbreak 

Therefore we have the following: 

\smallbreak   
$
X^{-1}=
\begin{pmatrix}
-1&-1\\
0&-1
\end{pmatrix}, 
$
$
^t\hskip-1mm X X^{-1}=
\begin{pmatrix}
1&1\\
-1&0
\end{pmatrix}, 
$
$
(^t\hskip-1mm X X^{-1})^2=
\begin{pmatrix}
0&1\\
-1&-1
\end{pmatrix}, 
$
$
(^t\hskip-1mm X X^{-1})^3=
\begin{pmatrix}
-1&0\\
0&-1
\end{pmatrix}.  
$

\bigbreak 
It is true that this knot exists. See e.g. \cite{Levinecob}. 
Furthermore we can say that 
this $(4k+1)$-knot $R$ 
is the knot product of 
the trefoil knot and 
$k$ copies of the Hopf link. 
See \cite{KauffmanOgasa} for the knot product.

\bigbreak
Let $C=(R\# (-R^*))\#(R\# (-R^*))$. 
Note that Arf $C=0$. 
By \cite{Ogasa98n}, 
$C$ is $(2k+1,2k+1)$-pass-move-equivalent to the trivial $(4k+1)$-knot. 

By Proposition \ref{Alabama} 
$H_{2k+1}(N_3(C);\Z)\cong\oplus^8\Z_2$. 
Hence we need no less than eight generators to generate $H_{2k+1}(N_3(C);\Z)$. 

Suppose that 
the $(2k+1,2k+1)$-pass-move-unknotting-number of $C$ is$\leqq1$.
By Proposition \ref{hsix} 
$H_{2k+1}(N_3(C);\Z)$ can take six generators. 
We arrived at a contradiction. 

Therefore 
the $(2k+1,2k+1)$-pass-move-unknotting-number of $C$ is$\geqq2$.

We prove that the $(2k+1,2k+1)$-pass-move-unknotting-number of $R\# (-R^*)$ is one. 
{\it Reason}: 
A $(2k+1)$-Seifert matrix of the $(4k+1)$-knot $R\# (-R^*)$ 
is the same as 
a 1-Seifert matrix of the uppermost 1-knot in Figure \ref{kirin}.  
(We have $R\# (-R^*)$ \newline
=(the uppermost 1-knot in Figure \ref{kirin})$\otimes^k$(the Hopf link), where $\otimes$ denotes the knot product which is defined 
in \cite{Kauffman, KauffmanNeumann}.)
One $(2k+1,2k+1)$-pass-move can change 
$R\# (-R^*)$ into a $(4k+1)$-knot $J_T$  
whose Seifert matrix is the same as 
that of the lower most 1-knot in Figure \ref{kirin}. 
Since the lower most 1-knot  in Figure \ref{kirin} is the trivial 1-knot, 
the Seifert matrix of $J_T$ is  $S$-equivalent to that of the trivial $(4k+1)$-knot 
(See \cite{Levinesimp} for  $S$-equivalence.)
By \cite{Levinesimp}, $J_T$ is the trivial  $(4k+1)$-knot. 
Hence  
the $(2k+1,2k+1)$-pass-move-unknotting-number of $R\# (-R^*)$ is one. 
(Another proof is given by using Main Theorem 2.4 of \cite{KauffmanOgasaII}.)

Therefore 
the $(2k+1,2k+1)$-pass-move-unknotting-number of $C$ is$\leqq2$.

Therefore 
the $(2k+1,2k+1)$-pass-move-unknotting-number of $C$ is two.

This completes the proof of Theorem \ref{Iowa}.(1). 
\bigbreak

Let $n\in\N$.  
Let $\#^nC$ be the connected-sum of $n$ copies of $C$. 
Since $C$ is $(2k+1,2k+1)$-pass-move-equivalent to the trivial $(4k+1)$-knot, 
$\#^nC$ is $(2k+1,2k+1)$-pass-move-equivalent to the trivial $(4k+1)$-knot. 

Let $N_3(\#^nC)$ be the 3-fold branched cyclic covering space of $S^{4k+3}$ along $\#^nC$.
By Proposition \ref{Alabama}, 
$H_{2k+1}(N_3(\#^nC);\Z)\cong\oplus^{8n}\Z_2$. 
We need no less than $8n$ generators which generate  $H_{2k+1}(N_3(\#^nC);\Z)$. 

Suppose that 
the $(2k+1,2k+1)$-pass-move-unknotting-number of $\#^nC$ is$\leqq n$. 
By Proposition \ref{hgsix}, we can prove that  
$H_{2k+1}(N_3(\#^nC);\Z)$ can take $6n$ generators.  
We arrived at a contradiction. 

Therefore 
the $(2k+1,2k+1)$-pass-move-unknotting-number of $\#^nC$ is$>n$.  

This completes the proof of Theorem \ref{Iowa}.(2). 

This completes the proof of Theorem \ref{Iowa}. \qed 

\begin{note}\label{penultimate}
Theorem \ref{Iowa} is a high-dimensional analogue of Theorem \ref{Arkansas}. 
Since  Theorem \ref{Arkansas} has a condition on the crossing-change-number, 
we naturally hope to impose a condition of a local-move which is a generalization of the crossing-change, on Theorem \ref{Iowa}. 
This is discussed in Note \ref{last}. 
\end{note}

\section{The twist-move on high-dimensional knots}\label{TM}   
\noindent 
Let $p\in\N\cup\{0\}$.  
We review the definition of 
the twist-move on $(2p+1)$-dimensional closed oriented submanifold$\subset S^{2p+3}$, 
which is defined in \cite{Ogasa09} 
and 
which is researched in \cite{KauffmanOgasa, KauffmanOgasaB, Ogasa09}.     
 (\cite{Ogasa09}  calls the twist-move the $XXII$-move.)  
If $p=0$,   
the twist-move on $(4p+1)$-dimensional closed oriented submanifold$\subset S^{2p+3}$
is the crossing-change on 1-links.

\begin{figure}
\vskip-5mm
\includegraphics[width=81mm]{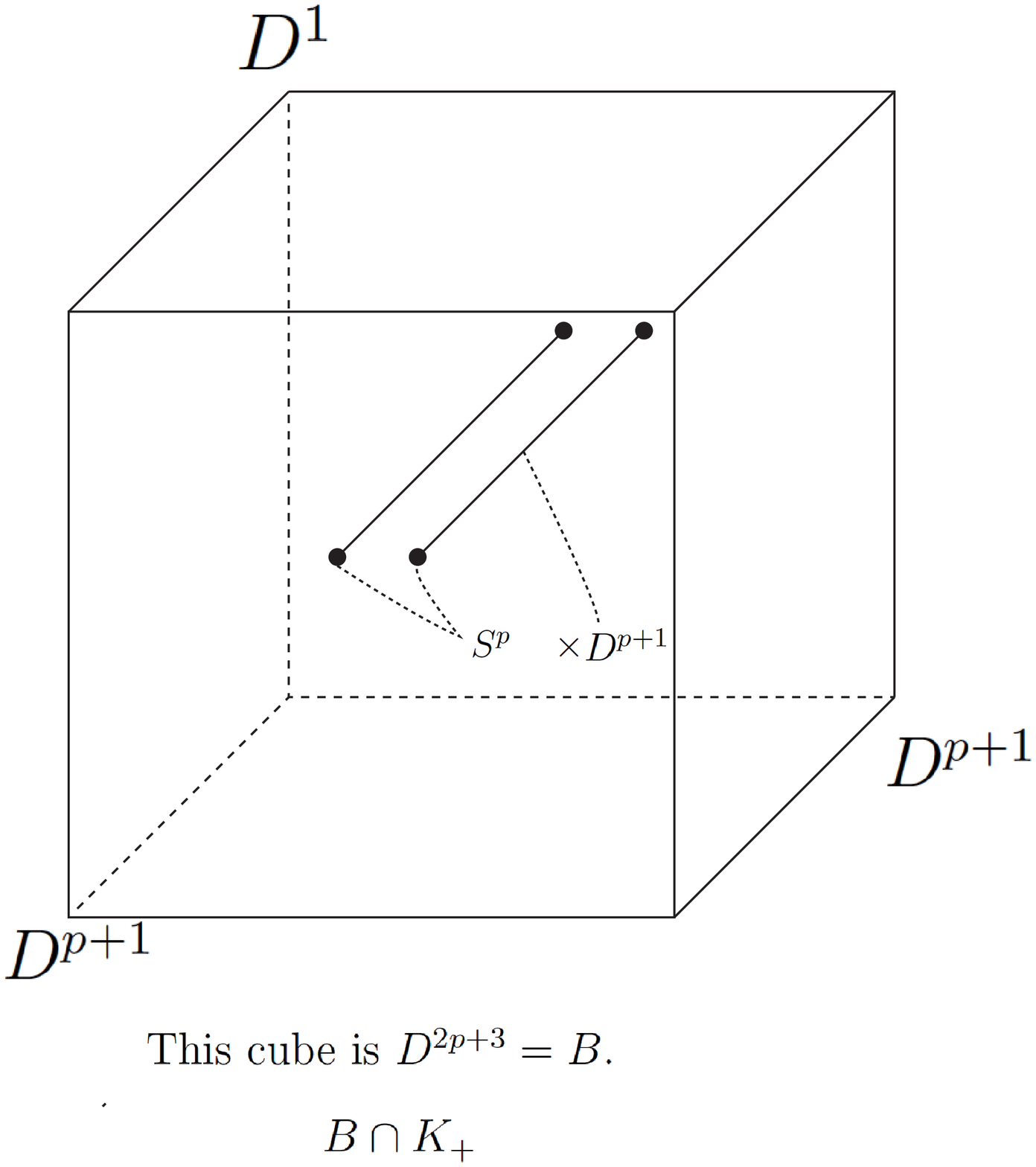}
\vskip-5mm
\caption{\hskip-2mm(1).{\bf The twist-move-triple}\label{chabin}}

\vskip12mm

\includegraphics[width=81mm]{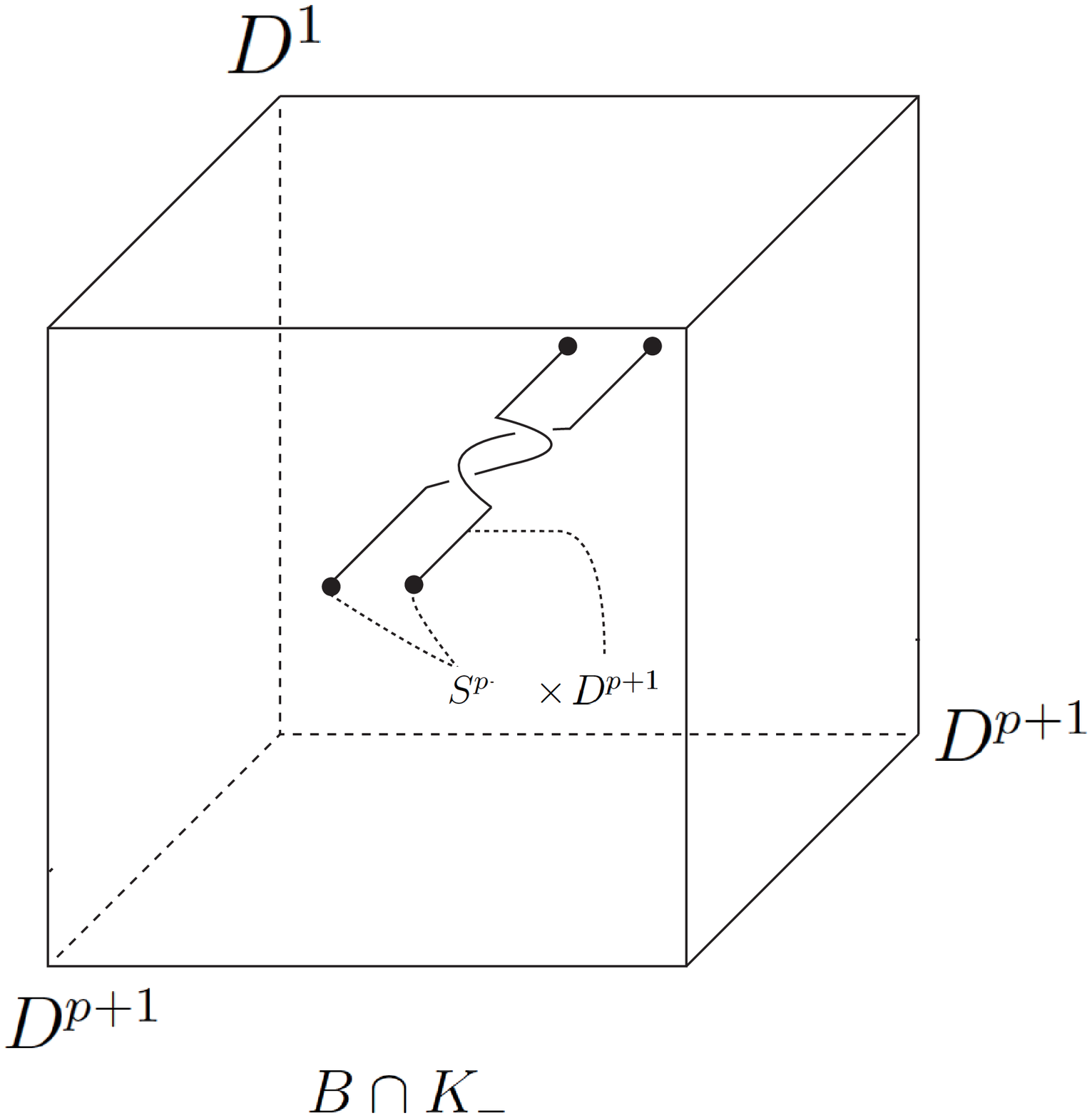}
\vskip-9mm
\text{F{\tiny IGURE} \ref{chabin}.(2): {\bf The twist-move-triple} }
\end{figure}

\begin{figure}
\includegraphics[width=11cm]{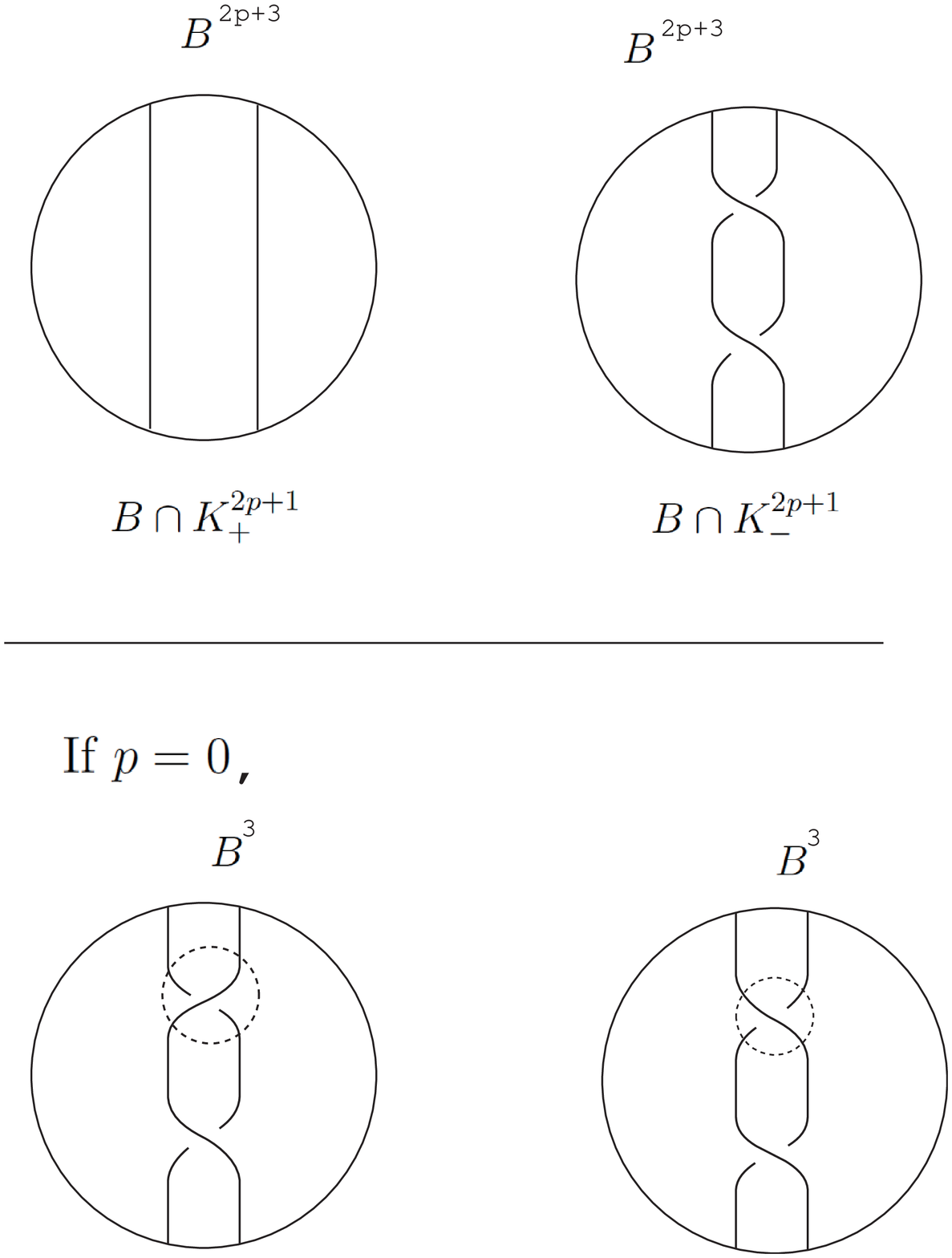}

\vskip2mm
\hskip20mm
 the pair of two   
\unitlength 0.1in
\begin{picture}(2.40,2.40)(0.60,-2.60)
%
\special{pn 8}%
\special{ar 180 140 120 120  0.0000000 0.1000000}%
\special{ar 180 140 120 120  0.4000000 0.5000000}%
\special{ar 180 140 120 120  0.8000000 0.9000000}%
\special{ar 180 140 120 120  1.2000000 1.3000000}%
\special{ar 180 140 120 120  1.6000000 1.7000000}%
\special{ar 180 140 120 120  2.0000000 2.1000000}%
\special{ar 180 140 120 120  2.4000000 2.5000000}%
\special{ar 180 140 120 120  2.8000000 2.9000000}%
\special{ar 180 140 120 120  3.2000000 3.3000000}%
\special{ar 180 140 120 120  3.6000000 3.7000000}%
\special{ar 180 140 120 120  4.0000000 4.1000000}%
\special{ar 180 140 120 120  4.4000000 4.5000000}%
\special{ar 180 140 120 120  4.8000000 4.9000000}%
\special{ar 180 140 120 120  5.2000000 5.3000000}%
\special{ar 180 140 120 120  5.6000000 5.7000000}%
\special{ar 180 140 120 120  6.0000000 6.1000000}%
\end{picture}%
   \hskip1mm
makes a crossing-change  on a 1-dimensional link.

\caption{{\bf  
The twist-move on 1-knots is the crossing-change-triple on 1-knots.
}\label{dobin}}
\end{figure}

\begin{defn}\label{spoon}
Let $p\in\N\cup\{0\}$.  
Regard a $(2p+3)$-ball $B=D^{2p+3}$ as $D^1\x D^{p+1}\x D^{p+1}$. 
Let $D^1=[-1,1]=\{t|-1\leqq t\leqq1\}$. 
Take a $p$-ball $D_S^{p+1}$ embedded in Int$D^{p+1}$. 
Let a submanifold $\{0\}\x D^{p+1}\x D_S^{p+1}\subset B$ 
be called $h_+$. 
We give an orientation to $h_+$. 
Take a submanifold $h_-\subset B$ which is diffeomorphic to   $h_+$. 
Let $h_+\cap h_-=h_+\cap(\partial B)$. 
Let $h_--(\partial B)\subset\{t<0\}\x D^{p+1}\x D^{p+1}$. 
(See 
Note (1) below.)
We give an orientation to $h_-$ so that 
 $h_+\cup h_-$ is an oriented submanifold$\subset B$  
if we give the opposite orientation to $h_-$.
We can regard $h_+\cup h_-$ as a Seifert hypersurface for $\partial(h_+\cup h_-)$. 
We can suppose that 
a $(p+1)$-Seifert matrix 
for a $(2p+1)$-dimensional closed oriented submanifold $\partial(h_+\cup h_-)\subset B$ 
associated with a Seifert hypersurface $h_+\cup h_-$ is $(1)$. 
(We can define Seifert hypersurfaces in $B$ and their Seifert matrices  
in the same fashion as ones in the $S^n$ case. 
Each of Figure \ref{chabin} and Figure \ref{dobin} draws a diagram of the twist-move.  
See 
Note (2) below.)

Let  $p\in\N\cup\{0\}$. 
Let $K_+$ and $K_-$ be 
$(2p+1)$-dimensional closed oriented submanifold $\subset S^{2p+3}$.
Take $B$ in $S^{2p+3}$. 
Let $K_+$ and $K_-$ differ only in $B$.
Let $K_+$ (resp. $K_-$) satisfy the condition  \newline
\hskip39mm$K_+\cap \mathrm{Int}B=
(\partial h_+-\partial B)$ \newline  
\hskip28mm 
${\rm(resp.}\hskip1mm K_-\cap \mathrm{Int}B=
(\partial h_--\partial B),$  \newline
where we suppose that there is not 
$h_-$  
(resp. $h_+$)   
in $B$.  
Then we say that 
$K_+$ (resp. $K_-$)  is obtained from $K_-$ (resp. $K_+$) 
by one 
{\it $($positive-$)$twist-move} 
(resp. {\it $($negative-$)$twist-move})  
in $B$.   

\end{defn}

In Definition \ref{spoon} we have the following: 
Let $\sharp\in\{+,-\}.$   
there is a Seifert hypersurface $V_\sharp\subset S^{2p+3}$ for $K_\sharp$ such that 
$V_\sharp\cap B=h_\sharp.$ 
\noindent 
(The idea of the proof is Thom-Pontrjagin construction.)
We say that  
$V_-$ (resp. $V_+$) is obtained from $V_+$ (resp. $V_-$) 
by one 
{\it $($positive-$)$twist-move} 
(resp. {\it $($negative-$)$twist-move})  
in $B$.

\bigbreak
\noindent
{\bf Note.} (1) \cite{Haefligerunknot, Haefligerknot, Whitney, Whitneytrick} etc. 
imply that 
we can move 
the core of  $h_-$ to 
the core of $h_+$ 
 in $B$ 
by an isotopy keeping $\partial$(the core of  $h_-$). 

\begin{figure}
\includegraphics[width=11cm]{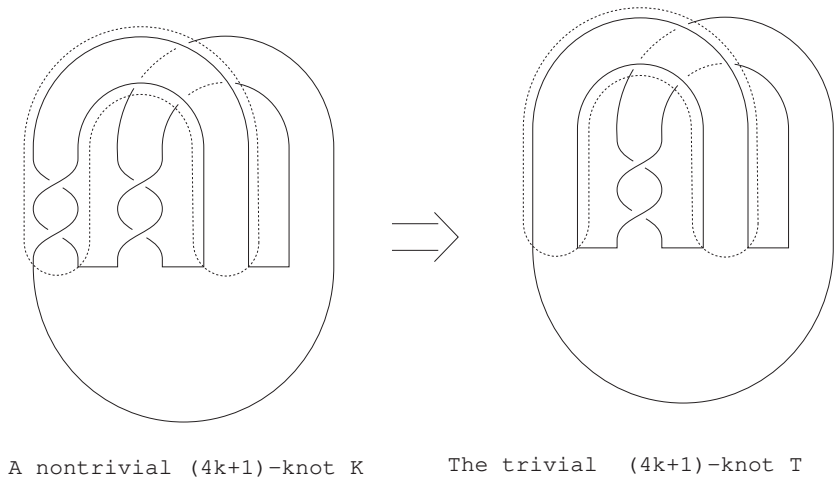}

\caption{{\bf 
One twist-move changes a nontrivial $(4k+1)$-knot 
into the trivial $(4k+1)$-knot
}\label{donburi}}
\end{figure}

\smallbreak
\noindent(2) 
Figure \ref{chabin}, which consists of the two figures (1) and (2),  
is a diagram of the twist-move. 
In Figure \ref{chabin}.(2), we move 
$\partial h_--\partial B$  
by isotopy and 
draw $\partial h_--\partial B$.  
The upper half 
of Figure \ref{dobin} is another diagram of the twist-move.
Compare 
the upper half 
of Figure \ref{dobin}  
and 
the lower half. 
If $p=0$ (hence $n=2p+1=1$), 
the left figure in the upper half
and 
that in the lower half are the same. 
That is, if $p=0$, 
the twist-move on $(2p+1)$-closed oriented submanifold$\subset S^{2p+3}$  
is the crossing-change on 1-links.  
Note that 
`$B\cap K_0$ in the left $B$ in the upper half 
of Figure \ref{dobin} in the $p=0$ case'  
and 
`$B\cap K_0$ in the left $B$ in the lower half of Figure \ref{dobin}'  
are the same ({\it Reason.} Use an isotopy.) 
See also Figure \ref{donburi}.  

\bigbreak
In Definition \ref{spoon}, note  the following: 
Let $\sharp\in\{+,-\}.$   
Let 
$V_0=V_\sharp-\text{Int}B$ \newline
$=\text{`the closure of }(V_\sharp- (h_\sharp))\text{ in }{S^{2p+3}}'.$  
We can say that 
we attach 
an embedded $(2p+2)$-dimensional $(p+1)$-handle $h_\#$ 
to the submanifold $V_0\subset S^{2p+3}$, 
and obtain the submanifold $V_\#\subset S^{2p+3}$.

\begin{defn}\label{atama}  
Let $p\in\N\cup\{0\}.$  
Let $K$ be a $(2p+1)$-knot$\subset S^{2p+3}$ 
which is twist-move-equivalent to the trivial $(2p+1)$-knot. 
The {\it twist-move-unknotting-number} of $K$ is 
the minimal number of twist-moves which we change $K$ to the trivial (2p+1)-knot by. 
\end{defn}




\begin{prop}\label{mimi}
There is a $(2p+1)$-knot whose twist-move-unknotting-number is one 
for a natural number $p$. 
\end{prop}


\noindent{\bf Proof of Proposition \ref{mimi}.} 
Let $k\in\N\cup\{0\}$. 
Take a $(4k+1)$-knot $K$ with a $(2k+1)$-Seifert matrix 
$
\begin{pmatrix}
-1&1\\
0&-1
\end{pmatrix}
$. See Figure \ref{donburi}.

Use the van Kampen theorem, the Mayor-Vietoris exact sequence, and \cite{Smale}. 
$K$ is PL homeomorphic to the standard sphere. 
(Note: Let $J$ be a $(2p+1)$-knot with a $(p+1)$-Seifert matrix 
$
\begin{pmatrix}
-1&1\\
0&-1
\end{pmatrix}
$ 
Then $J$ is not a homology sphere 
if $p$ is an odd natural number. {\it Reason.} Use the Mayor-Vietoris exact sequence.) 

By \cite{Levinesimp},  $K$ is a nontrivial spherical knot. 

We carry out one twist-move and 
obtain a $(4k+1)$-knot $T$ with a $(2k+1)$-Seifert matrix 
$
\begin{pmatrix}
0&1\\
0&-1
\end{pmatrix}
$.  
See Figure \ref{donburi}. 
By \cite{Levinesimp},  $T$ is the trivial knot. 

This completes the proof of Proposition \ref{mimi}.  
\qed

\bigbreak\noindent{\bf Note.}
By \cite{KervaireMilnor}, we have the following: 
$K$ in Figure \ref{donburi}  
is diffeomorphic to the standard sphere 
if  
$bP_{4k+2}$ is the trivial group.  
$K$ in Figure \ref{donburi}  
is diffeomorphic to an exotic sphere 
if  
$bP_{4k+2}$ is nontrivial.  
See  \cite{KervaireMilnor} for the ${\rm bP}$-subgroup.

\bigbreak
We consider the following problem.

\begin{prob}\label{kao}   
 Let $k\in\N\cup\{0\}.$   

\smallbreak\noindent
$(1)$ Is there a 
twist-move-unknotting-number-two $(4k+1)$-knot? 

\smallbreak\noindent
$(2)$ For any natural number $n$, 
is there a 
$(4k+1)$-knot 
whose twist-move-unknotting-number is$>n$?  
\end{prob}

We give a positive answer to Problem \ref{kao}.(1) (resp. \ref{kao}.(2)). 
The answers make one of our main theorems.   

\begin{thm}\label{kuchi}  
Let $k\in\N\cup\{0\}.$   Let $n\in\N$. 
%
There is a $(4k+1)$-knot whose twist-move-unknotting-number is $n$. 
\end{thm}

\bigbreak
\noindent{\bf Proof of Theorem \ref{kuchi}.} 
The $4k+1=1$ case holds because of 
\cite[Theorem 10.1 in page 420]{Murasugi}.  
Let $n\in\N$. 
The ordinary-unknotting-number, which is the twist-move-number,  of 
the connected-sum of $n$-copies of the trefoil knot
is $n$. 

The $4k+1\geqq5$ case is proved in the same fashion as 
one in 
\cite[Proof of Theorem 10.1  in page 420 and (2.4) in page 389]{Murasugi}.  
Let $k\in\N\cup\{0\}$. 
Take the $(4k+1)$-knot $K$ 
in Figure \ref{donburi} of Proof of Proposition \ref{mimi}.  
The twist-move-number of the connected-sum of $n$-copies $\sharp^nK$ of $K$ is $n$. 
\qed
\bigbreak
\noindent{\bf Note.}  
By \cite{KervaireMilnor}, we have the following: 
Let $n\in\N$. Let $k\in\N\cup\{0\}$. 
%
%
%
In the case where $bP_{4k+2}$ is nontrivial and $n$ is odd, 
$\sharp^nK$ is diffeomorphic to an exotic sphere.  
In the other case,  $\sharp^nK$ is diffeomorphic to the standard sphere.

\begin{note}\label{last}
We continue Note \ref{penultimate}. 
By the discussion in this section, 
we can prove that  
the twist-move-unknotting-number of $C$ in the previous section is $4$ 
and that that of $\#^nC$ is $4n$.  
It is natural to ask whether 
the twist-move-unknotting-number of  any $(4k+1)$-knot K is$\leqq4n$  
if the $(2k+1, 2k+1)$-pass-move-unknotting-number of $K$ is $n$. 
\end{note}

\bigbreak

\bigbreak\noindent
Eiji Ogasa:  Computer Science, Meijigakuin University, Yokohama, Kanagawa, 244-8539, Japan 
\quad pqr100pqr100@yahoo.co.jp  \quad
ogasa@mail1.meijigkakuin.ac.jp


\begin{thebibliography}{ABCD}

\bibitem{Browder} 
W. Browder: 
Surgery on simply-connected. manifolds, 
{\it Springer-Verlag Berlin Heidelberg New York.} (1972).











\bibitem{CochranOrr}  T. D. Cochran and K. E. Orr: 
Not all links are concordant to boundary links 
{\it Ann. of Math.} 138 (1993) 519--554. 





\bibitem{Haefligerunknot}  
A. Haefliger: Differentiable imbeddings, 
{\it Bull. Amer. Math. Soc.} 67 (1961) 109-112. 



\bibitem{Haefligerknot}  
A. Haefliger: Knotted $(4k - 1)$-spheres in $6k$-space, 
{\it Ann. of Math.} 75 (1962) 452-466. 






\bibitem{Kauffman}   
L. H. Kauffman:  Products of knots, 
{\it Bull. Amer. Math. Soc.} 
80 (1974)1104-1107.  

\bibitem{Kauffmanon}  
L. H. Kauffman: On knots, {\it Ann. of Math. Stud.} 115 (1987).

\bibitem{KauffmanNeumann}   
L. H. Kauffman and W. D. Neumann:  
Products of knots, branched fibrations and sums of singularities, 
{\it Topology} 16 (1977) 369-393.





\bibitem{KauffmanOgasa}   
L. H. Kauffman and E. Ogasa:  
Local moves of knots and products of knots, 
{\it 
Volume three of Knots in Poland III, 
Banach Center Publications} 103 (2014), 159-209, 
arXiv: 1210.4667 [math.GT]. 

\bibitem{KauffmanOgasaII}   
L. H. Kauffman and E. Ogasa:  
Local moves of knots and products of knots II,  
arXiv: 1406.5573 [math.GT].   


\bibitem{KauffmanOgasaB}   
L. H. Kauffman and E. Ogasa:  
Brieskorn submanifolds, Local moves on knots, and knot products, 
 arXiv: 1504.01229 [mathGT].  


\bibitem{KervaireMilnor}   M. Kervaire and J. Milnor: Groups of homotopy spheres I,  {\it Ann. of Math.} 77 (1963) 504-537.  



\bibitem{Kirby}
R. C. Kirby:  
The topology of 4-manifolds  {\it Springer-Verlag, Berlin, New York} (1989).







\bibitem{Levinepol}   
J. Levine: Polynomial invariant of knots of codimension two, 
{\it  Ann. of Math.} 84, (1966) 537-554.




\bibitem{Levinecob}  
J. Levine:  Knot cobordism in codimension two, 
{\it Comment. Math. Helv.} 44 (1969) 229-244. 


\bibitem{Levinesimp} J. Levine: 
An algebraic classification of some knots of codimension two. 
{\it Comment. Math. Helv.} 45 (1970) 185--198. 



\bibitem{Luck} 
W. L\"uck: A basic introduction to surgery theory. 
{\it ICTP Lecture Notes Series 9, Band 1, of the school "High-dimensional manifold theory" in Trieste, May/June $(2001)$, Abdus Salam International Centre for Theoretical Physics, Trieste.} 


\bibitem{Milnor}   
J.W.Milnor: Singular points of complex hypersurfaces, 
{\it Ann. of Math. Studies} 61. (1968).


\bibitem{Murasugi}
K. Murasugi: 
On a Certain Numerical Invariant of Link 
{\it Trans. Amer. Math. Soc.},  117, 387-422, 1965. 


\bibitem{Ogasa98n}   
E. Ogasa: 
Intersectional pairs of $n$-knots, local moves of $n$-knots and invariants of $n$-knots, 
{\it Math. Res. Lett.} 5 (1998) 577-582,  Univ. of Tokyo preprint UTMS 95-50.   

\bibitem{Ogasa02}     
E. Ogasa: The intersection of spheres in a sphere and 
a new geometric meaning of the Arf invariants, 
{\it J. Knot Theory Ramifications} 11 (2002) 1211-1231, 
Univ. of Tokyo preprint series UTMS 95-7, 
arXiv: 0003089 [math.GT].   



\bibitem{Ogasa04} 
E. Ogasa: Ribbon-moves of 2-links preserve the $\mu$-invariant of 2-links, 
{\it J. Knot Theory Ramifications}13 (2004) 669--687, 
UTMS 97-35, arXiv: 0004008 [math.GT].  



\bibitem{Ogasa07} 
E. Ogasa: Ribbon-moves of 2-knots: The Farber-Levine pairing and 
the Atiyah-Patodi-Singer-Casson-Gordon-Ruberman $\tilde\eta$ invariant of 2-knots, 
{\it Journal of Knot Theory and Its Ramifications} 16  (2007)  523-543, 
arXiv: 0004007 [math.GT], UTMS 00-22, arXiv: 0407164 [math.GT]. 


\bibitem{Ogasa09}  
E. Ogasa: Local move identities for the Alexander polynomials of 
high-dimensional knots and inertia groups,  
{\it J. Knot Theory Ramifications} 18 (2009) 531-545,    
 UTMS 97-63, arXiv: 0512168 [math.GT]. 




\bibitem{OgasaT3} 
E. Ogasa: A new obstruction for ribbon-moves of 2-knots: 2-knots fibred by the punctured 3-torus and 2-knots bounded by the Poincar\'e sphere,  
arXiv: 1003.2473 [math.GT].  



\bibitem{OgasaIH}
E. Ogasa: An introduction to high dimensional knots,  
arXiv: 1304.6053 [math.GT]. 


\bibitem{OgasaZ}  
E. Ogasa: Local-move-identities for the $\Z[t,t^{-1}]$-Alexander polynomials of 2-links, 
the alinking number, 
and 
high dimensional analogues, 
arXiv:1602.07775 [math.GT]. 


\bibitem{Ranichi}
A. Ranicki: 
Algebraic and geometric surgery
{\it Oxford Mathematical Monograph $($OUP$)$} 2002, 2003, 2014.

\bibitem{Smale} S. Smale:  Generalized Poincar\'e conjecture in dimensions greater than four,   {\it Ann. of Math.} 74  (1961) 391--406.


\bibitem{Wall}  C. T. C. Wall: 
 Surgery on compact manifolds, 
{\it   Academic Press, New York and London} (1970). 








\bibitem{Whitney}  
H. Whitney: Differentiable Manifolds, 
{\it Ann. of Math.} 37 (1936)  645-680. 


\bibitem{Whitneytrick}  
H. Whitney: The self-intersections of a smooth $n$-manifold in 2$n$-space
{\it  Ann. of Math.} 45 (1944) 220-246. 




\bibitem{Zeeman} E. Zeeman:  
Twisting spun knots, {\it Trans. Amer. Math. Soc.},   
115, 471-495, 1965. 



\end{thebibliography}
\end{document}